\documentclass[a4paper,10pt]{article}
\usepackage{graphicx}
\usepackage[colorlinks=true,linkcolor=blue,citecolor=blue]{hyperref}
\usepackage{color}

\usepackage{lipsum}
\usepackage[textsize=tiny]{todonotes}
\usepackage{amsmath, amssymb, amsthm}
\usepackage{mathtools}
\usepackage{dsfont}
\usepackage{nicefrac}
\usepackage{bbm}
\usepackage[inline]{enumitem}
\usepackage{xspace}
\usepackage[binary-units=true]{siunitx}
\usepackage{booktabs}
\usepackage{bookmark}
\usepackage[margin=3cm]{geometry}
\usepackage[affil-sl]{authblk}
\usepackage{thm-restate}
\usepackage{marvosym}
\usepackage{academicons}
\newcommand{\orcidicon}[1]{%
  \href{https://orcid.org/#1}{\hspace{1mm}\includegraphics[width=10pt]{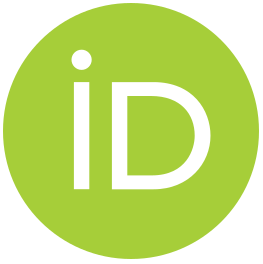}}
}

\usepackage{tikz}
\usetikzlibrary{arrows.meta}
\usetikzlibrary{fit, shapes.geometric, backgrounds}
\usetikzlibrary{decorations.pathreplacing, calligraphy}

\newcommand{\R}{\mathds{R}}
\newcommand{\Z}{\mathds{Z}}
\newcommand{\B}[1]{\{0,1\}^{#1}}

\newcommand{\solver}[1]{\textsc{#1}}
\newcommand{\scip}{\solver{SCIP}\xspace}
\newcommand{\setting}[1]{\texttt{#1}}

\newcommand{\fixed}[1]{\hat{ #1 }}
\newcommand{\lift}{l}

\newcommand{\sym}[1]{\mathcal{S}_{#1}}

\renewcommand{\log}[2][]{\text{log}_{#1}\left( #2 \right)}

\newcommand{\rectanglepack}{RP\xspace}
\newcommand{\spherepack}{SP\xspace}
\newcommand{\scheduling}{MS\xspace}

\newcommand{\sprod}[2]{#1^\top #2}
\newcommand{\define}{\coloneqq}
\newcommand{\card}[1]{\left\lvert #1 \right\rvert}

\newcommand{\brackets}[1]{\left\{ #1 \right\}}
\newcommand{\parentheses}[1]{\left( #1 \right)}

\newcommand{\naturalsto}[2][]{\left[ #2 \right]_{#1}}
\newcommand{\Oh}[1]{\mathcal{O}\parentheses{ #1 }}

\newtheorem{theorem}{Theorem}
\newtheorem{lemma}[theorem]{Lemma}
\newtheorem{remark}[theorem]{Remark}

\newtheorem{proposition}[theorem]{Proposition}

\theoremstyle{definition}
\newtheorem{definition}[theorem]{Definition}
\newtheorem{example}[theorem]{Example}
\newtheorem{observation}[theorem]{Observation}

\definecolor{col1}{HTML}{CD001A}
\definecolor{col2}{HTML}{EF6A00}
\definecolor{col3}{HTML}{79C300}
\definecolor{col4}{HTML}{F2CD00}
\definecolor{col5}{HTML}{1961AE}
\definecolor{col6}{HTML}{61007D}
\definecolor{col7}{HTML}{F27A7B}
\definecolor{col8}{HTML}{C27F01}
\definecolor{col9}{HTML}{3E882F}
\definecolor{col0}{HTML}{1F7231}

\begin{document}
\title{
  A Framework for 
  Handling and Exploiting Symmetry in Benders Decomposition\footnote{This article is part of the project ``Local Symmetries for Global Success'' with project number OCENW.M.21.299, which is financed by the Dutch Research Council (NWO).}
  }
\author{Christopher Hojny\orcidicon{0000-0002-5324-8996}}
\author{C\'edric Roy\orcidicon{0009-0007-4435-9508}\textsuperscript{(\Letter)}}
\affil{%
  Eindhoven University of Technology,
  Eindhoven, The Netherlands\\
  \emph{email} \{c.hojny, c.j.roy\}@tue.nl
}

\maketitle

%
%

\begin{abstract}
Benders decomposition (BD) is a framework for solving optimization problems by removing some variables and modeling their contribution to the original problem via so-called Benders cuts.
While many advanced optimization techniques can be applied in a BD framework, one central technique has not been applied systematically in BD yet: symmetry handling.
The main reason for this is that Benders cuts are not known explicitly but only generated via a separation oracle.

In this work, we close this gap by developing a theory of symmetry detection within the BD framework.
To this end, we introduce a tailored family of graphs that capture the symmetry information of both the Benders master problem and the Benders oracles.
Once symmetries of these graphs are known, which can be found by established techniques, classical symmetry handling approaches become available to accelerate BD.
We complement these approaches by devising techniques for the separation and aggregation of symmetric Benders cuts by means of tailored separation routines, extended formulations, and cut strengthening techniques.
All substantially reduce the number of executions of the separation oracles.
In a numerical study, we show the effect of symmetry handling, cut aggregation, and cut strengthening for bin packing and scheduling problems.
\end{abstract}

%
%

\section{Introduction}\label{sec:intro}

We consider mixed-integer programs (MIP)
\begin{equation}
  \min\{ \sprod{c}{x} + \sprod{d}{y} : Ax \leq b,\; Cx + Dy \leq f,\; x \in K(p,n),\; y \in K(q, m) \}, \label{eq:generic-MIP}
\end{equation}
where~$A$, $C$, $D$ are matrices and~$b$, $c$, $d$, $f$ are vectors of suitable dimensions, where~$K(r, s)$ is the shorthand notation for~$\Z^r \times \R^{s-r}$ for integers~$0\leq r\leq s$.
In many applications areas like healthcare delivery~\cite{heching2019logic}, sports scheduling~\cite{van2023traditional}, facility location~\cite{fazel2012using}, or packaging optimization~\cite{lehmann2023accelerated}, the matrix~$D$ has a block structure.
That is, when the~$x$-variables are fixed, the MIP decomposes into several independent smaller MIPs.

Benders Decomposition (BD)~\cite{benders1962} exploits this structure by splitting the MIP into a master and subproblem.
The master problem can be considered an abstract problem of the form
\begin{subequations}
  \label{eq:abstractBenders}
  \begin{align}
    \min \sprod{c}{x} + z & \label{eq:abstractBendersObjective}\\
    Ax &\leq b\label{eq:abstractBendersfluff}\\
    \varphi(x) &\leq z \label{eq:abstractBendersLink}\\
    x &\in K(p, n) \\
    z &\in \R
  \end{align}
\end{subequations}  
where the function~$\varphi(x) = \min\{\sprod{d}{y} : Dy \leq f - Cx,\; y \in K(q, m)\}$ models the optimal objective value of a~$y$-solution for a given~$x$-solution.
Variable~$z$ then models the contribution of the~$y$-variables to the objective and is linked to the~$x$-variables by~\eqref{eq:abstractBendersLink}.
Since~$\varphi$ does not admit a closed algebraic description, \eqref{eq:abstractBenders} cannot be solved by standard techniques.
Instead, \eqref{eq:abstractBendersLink} is linearized by a family of linear constraints~\eqref{eq:cutBendersLink}, where~$\mathcal{I}$ is some (usually finite) subset of~$\R^{n} \times \R$:
\begin{subequations}
  \label{eq:cutBenders}
  \begin{align}
    \min \sprod{c}{x} + z &&&\\
    Ax &\leq b&& \label{eq:bendersfluff}\\
    \sprod{\alpha}{x}  + \beta &\leq z, & (\alpha, \beta) &\in \mathcal{I} \label{eq:cutBendersLink} \\
    x &\in K(p, n)\\ 
    z &\in \R.
  \end{align}
\end{subequations}

To solve the original MIP, BD removes~\eqref{eq:cutBendersLink} from the master problem and adds, if needed, violated inequalities by calling a separation oracle.
The seminal paper~\cite{benders1962} by Benders showed how to implement a separation oracle for~\eqref{eq:cutBendersLink} when~$\varphi(x)$ is a linear program; logic-based BD extends these ideas to general MIPs~\cite{hooker2024,hooker2003logic}.
In particular, when~$D$ has a block structure, each block can have an independent oracle, which can reduce the running time of each oracle substantially. 

Once a separation oracle has been implemented, Problem~\eqref{eq:cutBenders} can be solved by standard branch-and-cut algorithms.
BD therefore also benefits from sophisticated MIP techniques (presolving, cutting planes, etc.), which are readily available in MIP software, with one important exception: symmetry handling.
The reason is that MIP software needs to have access to the full problem~\eqref{eq:cutBenders} to detect symmetries, whereas Inequalities~\eqref{eq:cutBendersLink} are only available through their  oracles.
Many researchers~\cite{adulyasakEtAl2015,akpinar2017combinatorial,delorme2017,guoEtAl2021} therefore add some symmetry handling inequalities by hand, but do not make use of advanced techniques such as orbital branching~\cite{OstrowskiEtAl2011} or orbitopal fixing~\cite{BendottiEtAl2021,DoornmalenHojny2024a,KaibelEtAl2011} to handle symmetries.
This motivates the following two questions:
\begin{enumerate*}[label={{(Q\arabic*)}},ref={Q\arabic*}]
\item\label{Q1} \emph{Can we develop a theory for detecting symmetries in BD?}
\item\label{Q2} \emph{Can we exploit symmetries in the generation of Benders cuts?}
\end{enumerate*}

\paragraph{Outline.}
To answer these questions, Section~\ref{sec:Benders-Symmetry} develops a practical framework for detecting symmetries in BD by using the recently introduced concept of symmetry detection graphs (SDG)~\cite{hojny2025detecting} to provide the master problem~\eqref{eq:cutBenders} symmetry information of the Benders oracles.
Detecting automorphisms of SDGs then allows to compute symmetries of~\eqref{eq:cutBenders} and to apply established symmetry handling methods.
In Section~\ref{sec:usingsymmetry}, we leverage symmetry information to enhance the separation of Benders cuts:
Let~$\pi$ be a symmetry of~\eqref{eq:cutBenders}.
If~$\sprod{\alpha}{x} + \beta \leq z$ is a Benders cut~\eqref{eq:cutBendersLink}, we observe that also~$\sprod{\pi^{-1}(\alpha)}{x} + \beta \leq z$ is a Benders cut.
A single Benders cut thus gives access to an entire family of symmetric Benders cuts.
For actions of the symmetric group, which arise in many applications, we show that this family of symmetric cuts can be separated very efficiently, thus avoiding to call potentially expensive separation oracles.
Moreover, by adding linearly many auxiliary variables and constraints we can compactly express all equivalent Benders cuts.
That is, separating symmetric cuts becomes obsolete once one member of the family is known.
In addition, we apply cut strengthening techniques to all the above-mentioned methods.
In Section~\ref{sec:Applications}, we illustrate how these techniques can be applied to bin packing and machine scheduling problems.
Section~\ref{sec:numerics} concludes the article with a numerical study in which we demonstrate that exploiting symmetry in BD can reduce the running time of BD by orders of magnitude.

This paper is an extended version of the IPCO conference paper~\cite{hojny2026IPCO}.
Next to more extensive numerical results, we present detailed mathematical proofs for the different techniques that were missing in the conference version.
Moreover, we extend the results by novel cut strengthening techniques that exploit symmetries in lifting procedures.

%
%

\section{Benders Symmetry}\label{sec:Benders-Symmetry}
In this section, we introduce our framework for detecting symmetries in BD.
We start by defining our notion of symmetries.
For a set~$A$, let~$\sym{A}$ be the set of all permutations of~$A$, the so-called \emph{symmetric group} of~$A$.
If~${A = [n] \define \{1,\dots,n\}}$ for some positive integer~$n$, we write~$\sym{n}$ instead of~$\sym{[n]}$.
For~$\pi \in \sym{n}$ and~${x \in \R^n}$, let $\pi(x) = (x_{\pi^{-1}(1)}, \dots, x_{\pi^{-1}(n)})$, i.e., $\pi$ acts on~$x$ by permuting its coordinates.
Following the MIP literature~\cite{liberti2012reformulations,margot2009symmetry}, a \emph{symmetry} of~$\min\brackets{\sprod{c}{x} : Ax \leq b, x\in K(p, n)}$ is a permutation~$\pi \in \sym{n}$ such that~$\sprod{c}{x} = \sprod{c}{\pi(x)}$ for all~$x \in \R^n$, and~$x$ is feasible if and only if~$\pi(x)$ is feasible. 
As deciding if~$\pi \in \sym{n}$ is a symmetry is NP-hard~\cite{margot2009symmetry}, one usually considers symmetries that keep a specific MIP formulation invariant.
Assuming~$A$ has~$m$ rows, a permutation~$\pi\in\sym{n}$ is a \emph{formulation symmetry}~\cite{Salvagnin2005} if there is~$\rho\in\sym{m}$ such that $\pi^{-1}(c) = c$, $\rho(b) = b$, and~$A_{\rho(i), \pi^{-1}(j)} = A_{i, j}$ for all~$(i, j)\in[n]\times[m]$.

To detect formulation symmetries, one usually translates the problem into the language of graph automorphisms.
To this end, one constructs a colored graph~$G = (V,E)$ with the following properties:
(i) there exists an injective map of the MIP's variables to a subset~$V'$ of~$V$, the \emph{variable nodes}; (ii) for every color-preserving automorphism~$\pi\colon V \to V$ of~$G$, the restriction of~$\pi$ to~$V'$ corresponds to a symmetry of the MIP.
We refer to such a graph as a \emph{symmetry detection graph (SDG)}.
For MIPs, \cite{Salvagnin2005} suggests a concrete construction of an SDG by introducing a node for each constraint~$i$ and variable~$j$.
These nodes are connected by an edge whose color corresponds to~$A_{i,j}$.
The constraint nodes~$i$ are assigned a color corresponding to their right-hand side~$b_i$.
Similarly, the variable nodes~$j$ are assigned a color encoding their objective coefficient~$c_j$ and their type (integer or continuous).
In the context of BD, however, this method is not applicable as the list of constraints~\eqref{eq:cutBendersLink} is often too large or not explicitly given.

Recently, \cite{hojny2025detecting} introduced a framework for constructing an SDG for a problem by combining SDGs of its building blocks.
A central concept for this construction are anchors.
An SDG~$G = (V,E)$ with variable nodes~$V'$ is called \emph{anchored}, if there exists a node~$a \in V \setminus V'$, the \emph{anchor}, such that, for every~$v \in V$, there exists an~$a$-$v$-path~$p$ in~$G$ such that no interior point of~$p$ belongs to~$V'$.
\begin{theorem}[cf. \cite{hojny2025detecting}]
  \label{thm:merge-SDGs}
  Let~$P = \min\{\sprod{c}{x} : Ax \leq b,\; Cx \leq d,\; x \in K(p,n)\}$ be a MIP.
  Let~$G_1$ and~$G_2$ be anchored SDGs of~$\min\{\sprod{c}{x} : Ax \leq b,\; x \in K(p,n)\}$ and~$\min\{\sprod{c}{x} : Cx \leq d,\; x \in K(p,n)\}$, respectively.
  Then, an SDG for~$P$ is given by the disjoint union of~$G_1$ and~$G_2$, and identifying the variable nodes of~$G_1$ and~$G_2$ with each other.
\end{theorem}
One can thus find an SDG for~\eqref{eq:abstractBenders} by merging anchored SDGs for~\eqref{eq:abstractBendersfluff} and~\eqref{eq:abstractBendersLink} (or~\eqref{eq:bendersfluff} and~\eqref{eq:cutBendersLink}).
This theorem forms the core of our framework for detecting symmetries in BD, which we illustrate for two different cases of BD.

\paragraph{Classical BD}
For illustration purposes, consider the linear program (LP)
\begin{subequations}
  \label{eq:trivial-example}
  \begin{alignat}{10}
    \min\       5&x_1  &+5&&x_2 &&+5&&y_1 &&+5y_2 &&+ 3y_3 &&+ 3y_4 \nonumber\\
                 &x_1  &+ &&x_2 &&  &&    &&      &&       &&       &&= 1,\\
                 &x_1  &  &&    &&+ &&y_1 &&      &&+ 2y_3 &&       &&= 2,\\
                 &     &  &&x_2 &&  &&    &&+ y_2 &&       &&+ 2y_4 &&= 2.
  \end{alignat}
\end{subequations}
If we keep only the~$x$-variables in the master problem, the master problem is given by
\begin{equation*}
  \min\{5x_1 + 5x_2 + z_1 + z_2 : x_1 + x_2 = 1,\; \varphi(x_1) \leq z_1,\; \varphi(x_2) \leq z_2\},
\end{equation*}
where~$\varphi(x_i) =\min\brackets{5y_i + 3y_{i+2} : y_i + 2 y_{i+2} = 2 - x_i}$ and~$i\in\brackets{1,2}$.
By Theorem~\ref{thm:merge-SDGs}, an SDG for the BD master problem is given by merging anchored SDGs for $\min\{5x_1+5x_2 + z_1 + z_2 : x_1 + x_2 = 1\}$ and~$\min\{5x_1 + 5x_2 + z_1 + z_2 : \varphi(x_i) \leq z_i\}$ for~$i \in \{1,2\}$, see Fig.~\ref{fig:SDG-mini-example} for an illustration.
Here, the dashed arcs correspond to the identification of variable nodes and anchors.
Since~\eqref{eq:cutBendersLink} is derived from~\eqref{eq:abstractBendersLink}, the SDG constructed this way is also an SDG for~\eqref{eq:cutBenders}.

This leaves the question how to find an SDG for \eqref{eq:abstractBendersLink} in general.
In the following, we show how to construct an SDG if the Benders subproblem is a linear program.
In our construction, we will also take into account a possible block decomposition that is exploited by BD.
To make this precise, let sets~$I_1, \dots, I_t$ partition~$[m]$.
If the matrix~$D$ in Problem~\eqref{eq:generic-MIP} is such that for any~$j\in I_{k}$ and~$j'\in I_{k'}$, $D_{jj'} = 0$ if~$k\neq k'$, then we can split the subproblem~$\varphi(x)$ into~$\varphi_k(x)$ for all blocks~$k\in [t]$, each only using variables~$\brackets{y_j: j\in I_k}$.
Problem~\eqref{eq:abstractBenders} can then be reformulated as
\begin{subequations}
  \label{eq:abstractBendersMultiBlock}
  \begin{align}
    \min \sprod{c}{x} &+ z_1 + \dots + z_t && \\
    Ax &\leq b,&&\\
    \varphi_k(x) &\leq z_k, && k \in [t], \label{eq:abstractBendersLinkBlock}\\
    x &\in K(p, n), && \\
    z &\in \R^t.&&
  \end{align}
\end{subequations}  
We can recover the symmetry detection graph of such formulations via the procedure described in the following statement.
\begin{theorem}
  \label{thm:SDG-Benders}
  Let~$G_1$ be an SDG for~\eqref{eq:generic-MIP} with variable nodes~$V^x\cup V^y_1 \cup \dots \cup V^y_k$, where the sets~$V^x = \brackets{v^x_1, \dots, v^x_n}$ and~$V^y_k = \brackets{v^y_j: j\in I_k}$ for any~$k\in[t]$ correspond to the variables~$x$ and~$y$ respectively.
  Let~$G_2$ be a modified copy of~$G_1$ with the following additions: 
  For each~$k\in [t]$ we add an uncolored variable node~$v^z_k$ connected to all nodes~$V^y_k$ with uncolored edges.
  Then $G_2$ is an SDG for the multi-block BD~\eqref{eq:abstractBendersMultiBlock}.
\end{theorem}
\begin{proof}
  Let~$\pi$ be a color-preserving automorphism of~$G_2$.
  We need to show that the restriction of~$\pi$ to the variable nodes is a symmetry of~\eqref{eq:abstractBendersMultiBlock}.
  To this end, let~$\bar{\pi}$ be the restriction of~$\pi$ to the nodes not in~$V^z = \{v^z_1,\dots,v^z_t\}$.
  Since~$V^z$ contains all uncolored variable nodes, $\pi$ keeps~$V^z$ invariant.
  Hence, $\bar{\pi}$ is an automorphism of~$G_2$ that only permutes nodes contained in~$G_1$.
  Therefore, $\bar{\pi}$ is an automorphism of~$G_1$.
  Because~$G_1$ is an SDG of~\eqref{eq:generic-MIP}, its restriction to the nodes corresponding to the $x$- and~$y$-variables is a symmetry of~\eqref{eq:generic-MIP}.

  To conclude the proof, we need to show that~$\varphi_k(x) = \varphi_{k'}(\pi(x))$ whenever~$\pi(v^z_k) = v^z_{k'}$ for some~$k,k' \in [t]$.
  If~$\pi(v^z_k) = v^z_{k'}$, then~$\pi(V^y_k) = V^y_{k'}$, because the nodes in~$V^y_k$ and~$V^y_{k'}$ are the only nodes adjacent with~$z^v_k$ and~$z^v_{k'}$, respectively.
  Here, we use the set notation~$\pi(S)$ for some set~$S$ to denote~$\{\pi(s) : s \in S\}$.
  As a consequence, $\varphi_k(x) = \varphi_{k'}(\pi(x))$, because the~$y$-variables defining~$\varphi_k$ are mapped onto those defining~$\varphi_{k'}$.
\end{proof}

\begin{figure}[tbp]
  \centering
  \begin{tikzpicture}[scale=0.75]
    \tikzstyle{cnb} += [draw=black, circle];
    \tikzstyle{cnw} += [draw=white, circle];
    \node (x1) at (  1, 1) [cnb, fill=col3, label=below:$x_1$]{};
    \node (x2) at (  2, 1) [cnb, fill=col3, label=below:$x_2$]{};
    \node (y1) at (  3, 1) [cnb, fill=col3, label=below:$y_1$]{};
    \node (y2) at (  5, 1) [cnb, fill=col3, label=below:$y_2$]{};
    \node (z1) at (  4, 1) [cnb, fill=col4, label=below:$y_3$]{};
    \node (z2) at (  6, 1) [cnb, fill=col4, label=below:$y_4$]{};
    \node (c1) at (  2, 4) [cnw, fill=col1]{};
    \node (c2) at (3.5, 4) [cnw, fill=col2]{};
    \node (c3) at (  5, 4) [cnw, fill=col2]{};
    \node (ak) at (3.5, 5) [cnw, fill=black!50]{};
    
    \begin{scope}[thick]
      \draw[col5] (x1) -- (c1);
      \draw[col5] (x2) -- (c1);
      \draw[col5] (x1) -- (c2);
      \draw[col5] (x2) -- (c3);
      \draw[col5] (y1) -- (c2);
      \draw[col5] (y2) -- (c3);
      \draw[col6] (z1) -- (c2);
      \draw[col6] (z2) -- (c3);  
      \draw[black!40] (c1) -- (ak);
      \draw[black!40] (c2) -- (ak);
      \draw[black!40] (c3) -- (ak);
    \end{scope}
  \end{tikzpicture}
  \hspace{1cm}
  \begin{tikzpicture}[scale=0.75]
    \tikzstyle{rnb} += [draw=black, rectangle];
    \tikzstyle{cnb} += [draw=black, circle];
    \tikzstyle{cnw} += [circle];
    \node (x_1) at (1,   1) [rnb, fill=col3, label=below:$x_1$]{};
    \node (x_2) at (2,   1) [rnb, fill=col3, label=below:$x_2$]{};
    \node (x_3) at (3, 3.5) [rnb, fill=col3, label=below:$x_1$] {};
    \node (x_4) at (6, 3.5) [rnb, fill=col3, label=below:$x_2$] {};
    \node (y_1) at (3, 2.5) [cnb, fill=col3, label=below:$y_1$]{};
    \node (y_2) at (6, 2.5) [cnb, fill=col3, label=below:$y_2$]{};
    \node (z_1) at (4, 2.5) [cnb, fill=col4, label=below:$y_3$]{};
    \node (z_2) at (7, 2.5) [cnb, fill=col4, label=below:$y_4$]{};
    \node (c_1) at (1.5, 4) [cnw, fill=col1]{};
    \node (c_2) at (  4, 4) [cnw, fill=col2]{};
    \node (c_3) at (  7, 4) [cnw, fill=col2]{};
    \node (w_1) at (  5, 1) [rnb, fill=cyan, label=below:$z_1$]{};
    \node (w_2) at (  8, 1) [rnb, fill=cyan, label=below:$z_2$]{};
    \node (ak1) at (2.5, 5) [draw=white, circle, fill=black!50]{};
    \node (ak2) at (  4, 5) [draw=white, circle, fill=black!50]{};
    \node (ak3) at (5.5, 5) [draw=white, circle, fill=black!50]{};

    \begin{scope}[thick]
      \draw[col5] (x_1) -- (c_1);
      \draw[col5] (x_2) -- (c_1);
      \draw[col5] (x_3) -- (c_2);
      \draw[col5] (x_4) -- (c_3);
      \draw[col5] (y_1) -- (c_2);
      \draw[col5] (y_2) -- (c_3);
      \draw[col6] (z_1) -- (c_2);
      \draw[col6] (z_2) -- (c_3);
      \draw[col0] (w_1) -- (c_2);
      \draw[col0] (w_2) -- (c_3);  
      \draw[black!40, dashed, -stealth] (x_4) to[out=225, in=0 ] (x_2);
      \draw[black!40, dashed, -stealth] (x_3) to[out=180, in=45] (x_1);
      \draw[black!40] (c_1) -- (ak1);
      \draw[black!40] (c_2) -- (ak2);
      \draw[black!40] (c_3) -- (ak3);
      \draw[black!40, dashed, -stealth] (ak2) to (ak1);
      \draw[black!40, dashed, -stealth] (ak3) to (ak2);
    \end{scope}
    
    \node (phi1) at (3.1, 4.4) [rotate=25]{$\varphi(x_1)$};
    \node (phi2) at (6.1, 4.4) [rotate=25]{$\varphi(x_2)$};
    \begin{scope}[on background layer]
        \filldraw[col7!30, line width=22pt] plot[smooth cycle] coordinates {(z_1.south) (y_1.south) (x_3) (c_2.south)};
        \filldraw[col7!30, line width=22pt] plot[smooth cycle] coordinates {(z_2.south) (y_2.south) (x_4) (c_3.south)};
    \end{scope}
  \end{tikzpicture}
  \caption{
    Anchored symmetry detection graphs for both formulations of~\eqref{eq:trivial-example}.
    On the left, SDG of the original formulation.
    On the right, example on how the SDG of the BD can be recovered. 
    We can join the SDGs of the subproblems~$\varphi(x_i)$, in the pink areas, and the SDG of the master problem by identifying the nodes~$x_i$.
    Variables of the master problem are depicted as squares nodes.
    }
  \label{fig:SDG-mini-example}
\end{figure}

\paragraph{Conflict Cuts}
Theorem~\ref{thm:merge-SDGs} can also be used to find an SDG for a Benders master problem if the subproblem is not an LP.
Note that we then refer to them as \emph{Logic Based Benders Decompositions} (LBBD)~\cite{hooker2003logic}.
These often appear, for example, in assignment problems, where we are given a set~$N$ of items and a set~$M$ of slots, and the task is to assign each item to exactly one slot while minimizing some cost~$c \in \R^{M \times N}$.
We can apply BD if the slots have the following monotonicity property: let~$S \subseteq S' \subseteq N$;
if assigning item set~$S\subseteq N$ to slot~$i$ is infeasible, then also the assignment of items~$S'$ to slot~$i$ is infeasible.
The Benders master problem with Benders cuts~\eqref{eq:cutBendersLink} can then be formulated as
\begin{subequations}
  \label{eq:naive-MK}
  \begin{align}
    \min \sprod{c}{x} & \\
    \sum_{i \in M} x_{ij} &= 1, && j \in N, \label{eq:multknapsack-partitioning}\\
    \sum_{j\in C} x_{ij} &\leq |C| -1, && (C, i) \in \mathcal{I}, \label{eq:multknapsack-BD} \\
    x &\in \B{M\times N}.
  \end{align}
\end{subequations}
The variable~$x_{ij}$ indicates whether item~$j$ is assigned to slot~$i$, and~\eqref{eq:multknapsack-partitioning} ensures that each item~$j$ is assigned to exactly one slot.
Constraint~\eqref{eq:multknapsack-BD} contains an inequality for every tuple~$(C,i)$, where~${C \subseteq N}$ is a set of items that cannot simultaneously be assigned to slot~$i \in M$.
We call these inequalities \emph{conflict cuts}, because they eliminate conflicting assignments of items.
We call the set~$P_i \define \brackets{(x_{i1}, \dots, x_{iN})\in\B{N} : x\text{ satisfies }\eqref{eq:multknapsack-BD}}$ the set of feasible \emph{partial assignments} to slot~$i$.
Depending on the application, the sets~$P_1, \dots, P_M$ have a simpler formulation.

\begin{example}
  \label{example:multipleknapsack}
  One example of the assignment problem is the multiple knapsack problem (MKP).
  Here, each item~$j \in N$ has a weight~$a_j$ and each slot~$i \in M$ has a capacity~$b_i$.
  Given a cost~$c_{ij}$ for assigning item~$j$ to slot~$i$,
  the task is to solve
  \begin{equation*}
    \min_{x\in \brackets{0, 1}^{M\times N}}\Big\{\sprod{c}{x} : 
    \sum_{i\in M} x_{ij} = 1 \text{ for all } j \in N \text{ and }
    \sum_{j\in N} a_j x_{ij} \leq b_i \text{ for all } i \in M\Big\}.
  \end{equation*}
  In a BD~\eqref{eq:naive-MK} of the MKP, one initially discards all conflict cuts, i.e., $\mathcal{I} = \emptyset$.
  The conflict cuts are then generated by repeatedly solving~\eqref{eq:naive-MK} and checking if the obtained solution~$\fixed{x}$ corresponds to a solution of the MKP.
  That is, the Benders oracles test if the sum of weights~$\sum_{j\in N}a_j\fixed{x}_{ij}$ exceeds the capacity~$b_i$ for some~${i \in M}$.
  In this case, the conflict cut~$(C,i)$ for~$C = \{j \in N: \fixed{x}_{ij} = 1\}$ is added to~$\mathcal{I}$; otherwise, the BD terminates and returns~$\fixed{x}$ as optimal~solution.
  As the Benders subproblem decomposes into problems~$P_i(x) = \min\{0 : \sum_{j \in N} a_j x_j \leq b_i\}$ for all~$i \in M$, an SDG for~\eqref{eq:naive-MK} is given by merging an SDG for~$\min\{\sprod{c}{x} : x \text{ satisfies~\eqref{eq:multknapsack-partitioning}}, x \in \B{M\times N}\}$ and SDGs for~$P_i(x)$, $i \in M$.
  The SDG for~$P_i(x)$ has a single constraint node connected to all variables nodes for~$x_{ij}$, $j \in N$, with edge weight~$a_j$, see Figure~\ref{fig:example-MPK}.
\end{example}

\begin{figure}[tb]
  \centering
  \begin{tikzpicture}[scale=0.8]
    \tikzstyle{cnb} += [draw=black, fill=black!70, circle, radius=5pt];
    \tikzstyle{cnw} += [draw=white, circle];
    \node (x11) at (  1, 1) [cnb, label=below:$x_{11}$]{};
    \node (x12) at (  2, 1) [cnb, label=below:$x_{12}$]{};
    \node (x13) at (  3, 1) [cnb, label=below:$x_{13}$]{};
    \node (x14) at (  4, 1) [cnb, label=below:$x_{14}$]{};
    \node (x15) at (  5, 1) [cnb, label=below:$x_{15}$]{};
    \node (x16) at (  6, 1) [cnb, label=below:$x_{16}$]{};
    \node (x21) at (  7, 1) [cnb, label=below:$x_{21}$]{};
    \node (x22) at (  8, 1) [cnb, label=below:$x_{22}$]{};
    \node (x23) at (  9, 1) [cnb, label=below:$x_{23}$]{};
    \node (x24) at ( 10, 1) [cnb, label=below:$x_{24}$]{};
    \node (x25) at ( 11, 1) [cnb, label=below:$x_{25}$]{};
    \node (x26) at ( 12, 1) [cnb, label=below:$x_{26}$]{};
    \node (x31) at ( 13, 1) [cnb, label=below:$x_{31}$]{};
    \node (x32) at ( 14, 1) [cnb, label=below:$x_{32}$]{};
    \node (x33) at ( 15, 1) [cnb, label=below:$x_{33}$]{};
    \node (x34) at ( 16, 1) [cnb, label=below:$x_{34}$]{};
    \node (x35) at ( 17, 1) [cnb, label=below:$x_{35}$]{};
    \node (x36) at ( 18, 1) [cnb, label=below:$x_{36}$]{};
    \node (P1)  at (3.5, 4) [cnw, fill=col3, label=above:$P_1$]{};
    \node (P2)  at (9.5, 4) [cnw, fill=col4, label=above:$P_2$]{};
    \node (P3)  at (15.5,4) [cnw, fill=col4, label=above:$P_3$]{};
    
    \begin{scope}[very thick]
      \draw[col5] (x11) -- (P1);
      \draw[col6] (x12) -- (P1);
      \draw[col6] (x13) -- (P1);
      \draw[col5] (x14) -- (P1);
      \draw[col1] (x15) -- (P1);
      \draw[col5] (x16) -- (P1);

      \draw[col5] (x31) -- (P3);
      \draw[col6] (x32) -- (P3);
      \draw[col6] (x33) -- (P3);
      \draw[col5] (x34) -- (P3);
      \draw[col1] (x35) -- (P3);
      \draw[col5] (x36) -- (P3);

      \draw[col5] (x21) -- (P2);
      \draw[col6] (x22) -- (P2);
      \draw[col6] (x23) -- (P2);
      \draw[col5] (x24) -- (P2);
      \draw[col1] (x25) -- (P2);
      \draw[col5] (x26) -- (P2);
    \end{scope}
  \end{tikzpicture}
  \label{fig:example-MPK}
  \caption{
    Example of SDG for MKP subproblems. 
    Every constraint node~$P_i$ has a color corresponding to its capacity~$b_i$.
    Every edge has a color corresponding to the weight~$a_i$ it is incident to, and all variable nodes get the same color.
    }
\end{figure}

\paragraph{Strengthening Conflict Cuts}
In general, conflict cuts are relatively weak, because they only cut off partial assignments related to a single item set~$C$.
In some cases, however, these cuts can be strengthened.
To illustrate this, consider a multiple knapsack problem in which all weights~$a_j$, $j \in N$, are identical.
In this case, once we know that~$(C, i)\in \mathcal{I}$ for some~$C\subseteq N$, we immediately get that~$(C',i) \in \mathcal{I}$ for every~$C' \subseteq N$ with~$|C| = |C'|$.
One can then deduce that the single cut~$\sum_{j\in N} x_{ij} \leq |C|- 1$ is valid and dominates the individual cuts for~$(C',i)$.
A standard approach to derive such a strengthened cut is lifting~\cite{padberg1975note}, which is frequently applied in integer programming.
Recently, we also showed how symmetries can be exploited for lifting knapsack based inequalities~\cite{hojny2020knapsack}.
But to the best of our knowledge, symmetry-aware lifting has not been exploited in the context of BD yet.
In the following, we provide the highlevel ideas of lifting.
The next section will then discuss how symmetries can be exploited in lifting conflict cuts.

For any given pair~$(C, i)\in\mathcal{I}$, one can \emph{lift}~\cite{kaparis2010separation} the corresponding Equation~\eqref{eq:multknapsack-BD} into a potentially stronger cut by determining valid coefficients for the variables not in the set~$C$.
In general, we want to find values~$\lift_{ij}$ for~$j\notin C$ such that
\begin{equation}
  \sum_{j\in C} x_{ij} + \sum_{j\notin C} \lift_{ij} x_{ij} \leq \card{C}-1
  \label{eq:LiftedNoGoodCut}
\end{equation}
is a valid cut.
One way to determine the value of these new coefficients is to solve an adjacent problem iteratively~\cite{padberg1975note,zemel1989easily}.
For any ordering of~$[n]\setminus C$, the lifting coefficient for~$j\notin C$ is computed via
\begin{equation}
  \label{eq:sequence-lifting}
  \lift_{ij} \define \min\brackets{
    |C| - 1 - \sum_{k\in C} x_{ik} - \sum_{k\in S} \lift_{ik} x_{ik} : x\in P_i,\ x_{ij} = 1
  },
\end{equation}
where~$S$ is the set of all elements coming before~$j$ in the chosen order.
The obvious downside to this technique is that it is sequence dependent.
Moreover, depending on the problem, the induced assignment instance can be very complicated and expensive to solve.
As we will see in the next section, however, exploiting symmetries can reduce the cost of lifting.

%
%

\section{Symmetry Exploitation}
\label{sec:usingsymmetry}
Once symmetries of~\eqref{eq:cutBenders} are known, a solver can make use of built-in symmetry handling methods.
Next to solving~\eqref{eq:cutBenders} by branch-and-bound, a computational bottleneck is the separation of Benders cuts~\eqref{eq:cutBendersLink}.
We therefore present two ways to exploit symmetries to generate Benders cuts.
The main idea revolves around the following observation.
\begin{observation}
  \label{obs:sym-cut}
  Let~$\pi \in \sym{n}$ be a symmetry of~\eqref{eq:cutBenders} and~$(\alpha,\beta) \in \mathcal{I}$ be a Benders cut of~\eqref{eq:cutBendersLink}.
  Then, $\sprod{\pi^{-1}(\alpha)}{x} + \beta \leq z$ is also a Benders cut.
\end{observation}
Indeed, for every feasible~$(x,z)$ of~\eqref{eq:cutBenders} and symmetry~$\pi$, we have~$z \geq \sprod{\alpha}{\pi(x)} + \beta = \sprod{\pi^{-1}(\alpha)}{x} + \beta$, where the last equality holds since permutations are orthogonal maps.
Every Benders cut~$(\alpha,\beta)$ thus gives rise to an entire family of symmetric Benders cuts~$\{(\pi^{-1}(\alpha),\beta) : \pi \text{ symmetry of~\eqref{eq:cutBenders}}\}$.
One can thus potentially avoid solving an expensive separation problem for~\eqref{eq:cutBendersLink} by first separating symmetric Benders cuts of already separated cuts.
In general, however, finding symmetric Benders cuts is intractable as deciding if there exists a~$\pi\in\Pi$ that separates a fixed~$\fixed{x},\fixed{z}$ is equivalent to deciding whether $\max_{\pi\in\Pi}\brackets{\sprod{\pi^{-1}(\alpha)}{\fixed{x}}} > \fixed{z} -\beta$.
\begin{proposition}[\cite{buchheim2005linear}]
  Let~$\Gamma \leq \sym{n}$, $a \in \R^{n}$, and~$x \in \R^n$.
  Then, computing~$\max_{\gamma \in \Gamma} \sprod{a}{\gamma(x)}$ is NP-hard.
\end{proposition}
In many applications, however, the symmetry groups of MIPs are well-structured, which allows to devise efficient methods for separating symmetric Benders cuts.
An important example is when the variables of an optimization problem can be arranged in an~$m \times n$ matrix~$X$ in such a way that the symmetries reorder the columns of~$X$.
That is, a column index permutation~$\gamma \in \sym{n}$ acts on~$X$ as~$(X_{i,j})_{i,j} \mapsto (X_{i,\gamma^{-1}(j)})_{i,j}$.
We refer to these symmetries as \emph{column symmetries}, which are the most frequently arising symmetries in practice, see~\cite{PfetschRehn2019}.
In fact, column symmetries naturally arise for assignment problems~\eqref{eq:naive-MK}, because the slots~$M$ are identical, i.e., we can arbitrarily exchange the sets of items assigned to the different slots.
Moreover, when fixing a set of symmetric items, the columns corresponding to these items can be permuted by column symmetries.
For such column symmetries, \cite{bendotti2026two} showed that symmetric cutting planes can be separated efficiently.
For the sake of completeness, we provide a proof sketch.
\begin{proposition}[\cite{bendotti2026two}]
  \label{prop:column-sym}
  Let~$A, X \in \R^{m \times n}$.
  Then, $\max_{\gamma \in \sym{n}}\sum_{i = 1}^m\sum_{j = 1}^n A_{i,j}X_{i,\gamma^{-1}(j)}$ can be computed in~$O(mn^2 + n^3)$ time.
\end{proposition}
\begin{proof}[Proof sketch.]
  Define the complete bipartite graph with node set~$\{u_1,\dots,u_n\} \cup \{v_1,\dots,v_n\}$ with edges connecting all~$u$-nodes with all~$v$-nodes.
  Assign edge~$\{u_j,v_{j'}\}$ the weight~$\sum_{i = 1}^m A_{i,j}X_{i,j'}$.
  The problem of maximizing~$\sum_{i = 1}^m\sum_{j = 1}^n A_{i,j}X_{i,\gamma^{-1}(j)}$ can then be solved by finding a weight maximum perfect matching in this graph.
  Constructing the weighted graph takes~$O(mn^2)$ time, whereas computing the matching can be done in~$O(n^3)$ time using the Hungarian method.
\end{proof}
As a consequence, symmetric Benders cuts for assignment problems~\eqref{eq:naive-MK} can be separated in polynomial time.
If there are many slots or identical items, however, the running time of this separation routine can still be rather high.
But by exploiting the structure of Benders cuts for~\eqref{eq:naive-MK}, we can substantially reduce the running time as we show next.
We note that, for some applications, this idea has been applied before, for example in network migration~\cite{Daryalal2023} or in railway resource scheduling~\cite{Han2025}, but it has not been described from a general perspective.
\begin{proposition}
  \label{prop:separateSn}
  Let~$\Pi$ be a symmetry group of~$P\define\min\brackets{\sprod{c}{x}: Ax\leq b, x\in K(p, n)}$.
  Let~$\alpha\in\R^n$, $\beta\in\R$ such that~$\sprod{\alpha}{x}\leq \beta$ is a valid cut of~$P$. 
  If~$\Pi = \sym{n}$, then~$\brackets{\sprod{\pi^{-1}(\alpha)}{\fixed{x}} \leq \beta : \pi\in\Pi}$ can be separated in~$\Oh{n\log{n}}$ time for any fixed~$\fixed{x}$.
\end{proposition}
\begin{proof}
  Let~$\fixed{x} \in K(p,n)$.
  Determining if a~$\pi\in\Pi$ separates~$\fixed{x}$ is equivalent to deciding whether $\max_{\pi\in\Pi}\sprod{\pi^{-1}(\alpha)}{\fixed{x}} > \beta$.
  Since~$\Pi=\sym{n}$, finding a maximizer~$\pi$ is the same as finding a permutation that matches the $k$-th largest value of~$\brackets{\alpha_j : j\in [n]}$ to the $k$-th largest value of~$\brackets{\fixed{x}_j : j\in [n]}$, which can be done by sorting both~$\alpha$ and~$\fixed{x}$ in~$\Oh{n\log[]{n}}$ time.  
\end{proof}

\begin{proposition}
  \label{prop:separate-benders}
  Let~$\Pi$ be a symmetry group of~\eqref{eq:naive-MK} and~$(\alpha, \beta) \in \mathcal{I}$ a Benders cut.
  If~$\Pi = \sym{n}$, then~$\{(\pi^{-1}(\alpha),\beta) : \pi\in\Pi\} \subseteq \mathcal{I}$ can be separated in~$\Oh{n\log{n}}$ time.
\end{proposition}
\begin{proof}
  Observe that although~$\Pi$ acts on the columns of~$x$, the cut~$(\alpha, \beta)\in\mathcal{I}$ only has non-zero values on a single row of~$x$.
  In other words, we can ignore the action of~$\Pi$ on the other rows of~$x$ and apply the same argument as in Proposition~\ref{prop:separateSn}.
\end{proof}
Note that the previous three propositions can be expanded to products of symmetric groups.
Let~$\mathcal{G}$ be a partition of~$[n]$.
If~$\Pi = \bigotimes_{G\in \mathcal{G}} \sym{G}$, i.e., $\Pi$ is a direct product of symmetric groups, we can apply the same proof on each~$G\in\mathcal{G}$ independently.
This allows us to consider instances where not all items are identical, but have rather groups of identical items.
For the remainder of this section, we will assume that the symmetry group~$\Pi$ of the problem at hand is a direct product of symmetric groups.

We refer to~$\{(\pi^{-1}(\alpha),\beta) : \pi\in\Pi\}$ as a \emph{symmetric cut pool} derived from~$(\alpha,\beta)$.
Figure~\ref{fig:cut-pool-scheme} illustrates how symmetric cut pools can be used in Benders decomposition.
After the master problem has been solved and a solution~$\fixed{x}$ has been obtained, we can cheaply try to separate cuts from the symmetric cut pool.
If~$\fixed{x}$ can be separated, we immediately add the corresponding cut and update the master problem.
Otherwise, the potentially more expensive separation method via the Benders subproblem needs to be solved to decide feasibility and optimality of~$\hat{x}$. 

\begin{figure}[tb]
  \centering
  \begin{tikzpicture}[scale=1]
    \node (master) at (0,2) [rectangle,draw=black] {solve master problem};
    \node (pool) at (4, 2) [circle, draw=black] {
      \begin{minipage}{0.05\linewidth}
        \begin{center} cut \\ pool? \end{center}
      \end{minipage}
    };
    \node (sepa) at (8,0) [rectangle,draw=black] {separate master solution};
    \node (viol) at (8,-2) [circle,draw=black] {
      \begin{minipage}{0.05\linewidth}
        \begin{center}found\\ cut?\end{center}
      \end{minipage}
      };
    \node (solved) at (10,-3) [rectangle,draw=black] {problem solved};
    \node (ac) at (4,-2) [rectangle, draw=black] {add cut};
    \node (viol2) at (4,0) [circle,draw=black] {
      \begin{minipage}{0.05\linewidth}
        found\\ sym.\\ cut?
      \end{minipage}
      };
    \node (aux) at (3,-1.5) [inner sep=0mm,minimum size=0mm] {};

    \begin{scope}[-stealth, thick]
      \draw (master) -- (pool);
      \draw (pool) -| node[pos=0.1, above] {no} (sepa);
      \draw (pool) -- node[pos=0.5, right] {yes} (viol2);
      \draw (viol2) -- node[pos=0.5, above] {no} (sepa);
      \draw (viol2) -- node[pos=0.5, right] {yes} (ac);
      \draw (sepa) -- (viol);
      \draw (viol) -| node[pos=0.25, above] {no} (solved);
      \draw (viol) -- node[pos=0.25, above] {yes} (ac); 
      \draw (ac) -| (master);
    \end{scope}
  \end{tikzpicture}
  \label{fig:cut-pool-scheme}
  \caption{
   Schematic illustration of the solving process for Benders decomposition.
   Observe that the use of the cut pool technique allows to bypassing the separation of the master solution via the subproblem for any solution of a known symmetric family.
   }
\end{figure}

Propositions~\ref{prop:column-sym}--\ref{prop:separate-benders} allow to possibly accelerate the separation of~\eqref{eq:cutBendersLink}.
A natural question is if the separation of symmetric Benders cuts can be avoided at all.
In the case of our assignment problems, this question can be answered affirmatively.
Extending ideas of~\cite{riise2016recursive}, we introduce for each group~$G \in \mathcal{G}$ variables~$\zeta^k_{G}\in \brackets{0, 1}$ and constraints that enforce~$\zeta^k_{G} = 1$ if and only if~$\sum_{j\in G} x_{j} \geq k$, enforced via the two following constraint sets:
\begin{align}
  \sum_{j \in G} x_j &\leq k - 1 + \left(|G| - (k - 1)\right)\zeta^k_G, && G \in \mathcal{G},\; k \in [|G|],\label{eq:zeta-lowerbound}\\
  k \cdot \zeta^k_G &\leq \sum_{j \in G} x_j, && G \in \mathcal{G},\; k \in [|G|].\label{eq:zeta-upperbound}
\end{align}
The variables~$\zeta^1_{G}, \dots, \zeta^{|G|}_{G}$ can thus be considered a non-increasing reordering of~$\brackets{x_{j}: j \in G}$.
By denoting the~$k$-th largest value in~$\{\alpha_j : j \in G\}$ by~$\alpha^k_{G}$, the most violated cut constructed in the proof of Proposition~\ref{prop:separate-benders} can be phrased in terms of the~$\zeta$-variables as
\begin{equation}
  \sum_{G\in \mathcal{G}} \sum_{k = 1}^{|G|} \alpha^k_{G} \zeta^k_{G} +\beta \leq z. \label{eq:generic-EF-cut}
\end{equation}
We thus immediately obtain the following theorem.
\begin{theorem}
  \label{thm:EF-works}
  Assume all~$x$-variables in~\eqref{eq:cutBenders} are binary.
  Let~$(\alpha,\beta) \in \mathcal{I}$ and symmetry group~$\Pi$ adhere to the conditions in Proposition~\ref{prop:separate-benders}.
  Then, $(\fixed{x},\fixed{z}) \in \B{n} \times \R$ satisfies all equivalent inequalities in~$\{(\pi^{-1}(\alpha),\beta) : \pi \in \Pi\}$ if and only if~$(\fixed{x},\fixed{\zeta},\fixed{z})$ in the extended model satisfies~\eqref{eq:generic-EF-cut}.
\end{theorem}

We now have two approaches to benefit from symmetries when generating Benders cuts: the symmetric cut pool and the extended formulation.
It is unclear though which of these two formulations is stronger.
In particular, the extended model is only an extended IP formulation, but we have not settled whether it is also an extended formulation of the LP relaxation of the cut pool model.
In the following, we answer this question rigorously.
While Theorem~\ref{thm:EF-works} shows that the separation of symmetric Benders cuts can be avoided, the extended model~\eqref{eq:zeta-lowerbound}--\eqref{eq:generic-EF-cut} has a weaker LP relaxation than~\eqref{eq:cutBenders} in general.
Indeed, Proposition~\ref{lemma:EF-is-weaker} shows that every feasible~$(x, z)\in [0,1]^{n}\times\R$ induces a feasible~$(x, \zeta, z)\in[0,1]^{2n}\times\R$.
The converse is however not true.
The next example shows that there exist some~$(x,\zeta, z)$ with projection to~$(x,z)$ not satisfying~$\sprod{\alpha}{x} + \beta \leq z$.

\begin{example}
  Consider a problem with three variables $x$, $y_1$, and $y_2$ for which the conflict cuts $x + y_1 \leq 1$ and $x + y_2 \leq 1$ can be derived.
  Assuming that $y_1$ and $y_2$ are symmetric, the two conflict cuts belong to a symmetric cut pool, which can be represented by $\zeta^1_x + \zeta^1_y \leq 1$ in the extended model.
  Then, for every $\varepsilon \in [0,\frac{1}{2}]$, the solution $(x,y_1,y_2,\zeta^1_x,\zeta^1_y,\zeta^2_y) = (1-\varepsilon, 2\varepsilon, 0, 1-\varepsilon,\varepsilon,\varepsilon)$ is feasible for the LP relaxation of the extended model.
  However, solution $(x,y_1,y_2) = (1-\varepsilon,2\varepsilon,0)$ is not feasible for the cut pool model whenever~$\varepsilon > 0$.
\end{example}

\begin{restatable}{proposition}{lemmaEFisWeaker}
  \label{lemma:EF-is-weaker}
  Let~$(x,z)\in[0, 1]^{n}\times\R$, $\alpha\in\R^n_{\geq 0}$ satisfy~$\sum_{j\in C} \alpha_j \pi(x)_j +\beta \leq z$ for all~$\pi\in\Pi$.
  Then there exists a~$\zeta$ such that~$(x,\zeta, z)$ is feasible in the relaxation of the extended formulation.
\end{restatable}

The proof of Proposition~\ref{lemma:EF-is-weaker} requires two intermediary results.
Namely, we first show that we can construct~$\zeta$ such that its entries are upper-bounded by the average of~$x$ within the corresponding group~$G$ (Lemma~\ref{rmk:zeta-lb}).
We can then show that such a~$\zeta$ will satisfy constraint~\eqref{eq:generic-EF-cut}.
Trivially, the sum of the first~$k$ largest element of~$x$ is always greater than~$k$ times the average.
We show that this property still holds in a weighted sum when the weights are non-increasing (Lemma~\ref{rmk:weighted-average-dom}).
\begin{lemma}
  \label{rmk:zeta-lb}
  Let~$x\in[0, 1]^n$. 
  There exists~$\zeta\in [0, 1]^n$ satisfying~\eqref{eq:zeta-lowerbound} and~\eqref{eq:zeta-upperbound} such that~$\zeta^k_G \leq \frac{x(G)}{|G|}$ for all~$G\in\mathcal{G}$ and~$k\leq|G|$.
\end{lemma}
\begin{proof}
  For ease of notation, let~$x(S)$ denote~$\sum_{j\in S} x_j$ for any set~$S\subseteq\naturalsto[]{n}$.
  Combining Equations~\eqref{eq:zeta-lowerbound} and~\eqref{eq:zeta-upperbound}, we get that
  \begin{equation*}
    \frac{x(G) -k + 1}{|G| - k + 1} \leq \zeta^k_G \leq \frac{x(G)}{k}
  \end{equation*}
  which is a non-empty interval for any~$1\leq k\leq |G|$ and~$x\in [0, 1]^n$.
  Indeed,
  \begin{equation*}
    \frac{x(G)}{k} - \frac{x(G) - k + 1}{|G|-k+1} = \frac{x(G)(|G| -2k+1) + k(k-1)}{k(|G|-k+1)}.
  \end{equation*}
  Observe that the numerator can be rewritten as a function of~$k$ or of~$x(G)$ in the following:
  \begin{equation*}
    x(G)\left(|G| + 1 - 2k\right) + k(k-1) = k^2 - (2x(G)+1)k + x(G)(|G|+1).
  \end{equation*}
  The minimum with respect to~$x(G)$ is either at~$x(G) = 0$ or~$|G|$, depending on whether~$k \geq \frac{|G|+1}{2}$. 
  Reciprocally, the minimum with respect to $k$ is when~$2k - (2x(G) + 1) = 0$, or equivalently when~$k$ is close to~$\frac{1}{2} + x(G)$.
  As~$k\in \{1, \dots, |G|\}$, the global minimum of the numerator is either when~$k = 1, x(G) = 0$ or~$k=x(G)=|G|$.
  In either case, the numerator becomes~$0$, meaning that the interval for~$\zeta^k_G$ is never empty.
  Let~$\hat{\zeta}^k_G$, $k\leq |G|$ be the unique~$\zeta$-values that satisfy~\eqref{eq:zeta-lowerbound} with equalities.
  We can observe that for all~$G\in\mathcal{G}$, $\hat{\zeta}^1_G = \frac{x(G)}{|G|}$ and for all~$k\leq|G|-1$ we have
  \begin{align*}
    \hat{\zeta}^k_G - \hat{\zeta}^{k+1}_G = \frac{x(G) - k + 1}{|G| - k + 1} - \frac{x(G) - k}{|G| - k} = \frac{|G| - x(G)}{(|G|-k+1)(|G|-k)} \geq 0.
  \end{align*}
  Observe that~$\hat{\zeta}^k_G$ satisfying~\eqref{eq:zeta-lowerbound} with equality does not always guarantee~$\hat{\zeta}^k_G \geq 0$.
  Nevertheless, as the upper bound derived from~\eqref{eq:zeta-upperbound} is always greater or equal to zero,  $\zeta^k_G = \max\brackets{\hat{\zeta}^k_G, 0}$ ensures that~$\frac{x(G)}{|G|} \geq \zeta^1_G \geq \zeta^2_G \geq \dots \geq \zeta^{|G|}_G$ as claimed.
\end{proof}

\begin{lemma}
  \label{rmk:weighted-average-dom}
  Let~$\alpha, x\in \R^n$ such that~$\alpha_1 \geq \dots \geq \alpha_{n} \geq 0$ and~$x_1 \geq \dots \geq x_{n} \geq 0$.
  Denoting the average of~$x$ by~$A_v(x)\define \sum^n_{k=1} x_i/n$, then $\sum_{k=1}^K \alpha_k A_v(x) \leq \sum_{k=1}^K \alpha_k x_k$ for any~$K\leq n$.
\end{lemma}
\begin{proof}
  Let~$f(K)\coloneq \sum_{k=1}^K \alpha_k (x_k - A_v(x))$. 
  We prove that~$f(K)\geq 0$ for any~$K\leq n$.
  Observe that~$f$ behaves like a concave function as~$\alpha_k$ and~$(x_k - A_v(x))$ both decrease when~$k$ increases and~$\alpha_k \geq 0$.
  We then only need to check that~$f(1) \geq 0$ and~$f(n)\geq 0$.
  Since~$x_1 \geq A_v(x)$ and $\alpha_1 \geq 0$, we conclude $f(1) \geq 0$.
  Let~$g(k)\coloneqq \sum_{i=1}^k (x_i - A_v(x))$.
  It is easy to verify that~$g(k)\geq 0$ for the same arguments:~$g$ also behaves like a  concave function, $g(1)\geq 0$, and~$g(n) = 0$.
  Notice that since~$g(k) - g(k-1) = x_k - A_v(x)$ for all~$k\in \brackets{2, ..., n}$, we have that
  \begin{align*}
    f(n)
    &= \sum_{k = 1}^n \alpha_k(x_k - A_v(x)) = \sum_{k = 2}^n \alpha_k(x_k - A_v(x)) + \alpha_1 (x_1 - A_v(x)) \\
    &= \sum_{k=2}^{n} \alpha_k (g(k) - g(k-1)) + \alpha_1 g(1) = \sum_{k=2}^{n} \alpha_k g(k) - \sum_{k=2}^{n} \alpha_k g(k-1) + \alpha_1 g(1) \\
    &= \alpha_n g(n) + \sum_{k=2}^{n-1} \alpha_k g(k) - \sum_{k=1}^{n-1} \alpha_{(k+1)} g(k) + \alpha_1 g(1) \\
    &= \alpha_n \underbrace{g(n)}_{=0} + \sum_{k = 1}^{n-1} (\underbrace{\alpha_k - \alpha_{k+1}}_{\geq 0}) \underbrace{g(k)}_{\geq 0} 
  \end{align*}
  and thus~$f(n) \geq 0$.
\end{proof}

We now have all ingredients to prove Proposition~\ref{lemma:EF-is-weaker}.
\begin{proof}[Proof of Proposition~\ref{lemma:EF-is-weaker}]
  Let~$x\in [0,1]^n$ with and $z \in \R$ satisfy~$\sum_{j\in C} \alpha_j \pi(x)_j \leq z - \beta$ for all~$\pi\in\Pi$.
  Recall that by the definition of~$\Pi$, this means that~$\sum_{G\in \mathcal{G}} \sum_{j\in C\cap G} \alpha_j \pi_G(x)_j \leq z - \beta$ for any permutations~$\pi_G\in\sym{G}$ and~$G\in\mathcal{G}$.
  Combining Lemma~\ref{rmk:zeta-lb} and~\ref{rmk:weighted-average-dom}, there exists a~$\zeta$ such that
  \begin{align*}
    \sum_{G\in\mathcal{G}} \sum_{k=1}^{|G\cap C|} \alpha^k_G \zeta^k_G 
    &\overset{\text{Lemma }\ref{rmk:zeta-lb}}{\leq} \sum_{G\in\mathcal{G}} \sum_{k=1}^{|G\cap C|} \alpha^k_G \frac{x(G)}{|G|} \\
    &\overset{\text{Lemma }\ref{rmk:weighted-average-dom}}{\leq} \sum_{G\in\mathcal{G}} \max_{\pi_G\in\sym{G}} \sum_{j\in C\cap G} \alpha_j \pi_G(x)_j 
    = \max_{\pi\in\Pi} \sum_{j\in C} \alpha_j \pi(x)_j \\
    &\overset{\phantom{\text{Lemma x}}}{\leq} z - \beta.
  \end{align*}
\end{proof}

Note that another, arguably more natural extended formulation could have used integer variables~$\xi_G\in\Z_{\geq 0}$ with~$\xi_G = \sum_{j\in G} x_j$.
It is however unclear how they can be applied in a cut.
In fact, since the convex hull of all feasible~$(x, \xi)$ can contain infeasible integer points, as illustrated in the following example, they would likely require additional auxiliary variables which would then be equivalent to the~$\zeta$ variable.
\begin{example}
  \label{example:xi-not-convex}
  Consider an assignment problem where one needs to fit non-overlapping rectangular boxes inside rectangular bins.
  Assuming the bins are of dimension~$K\times K$ and that the boxes are either of dimension~$K\times 1$ or~$1\times K$ for some~$K \in\Z_{>1}$.
  It is then obvious that a solution can involve up to~$K$ boxes of same dimension at the same time, but one cannot have a horizontal and vertical box simultaneously.
  In terms of the~$\xi = (\xi_h, \xi_v)\in\Z^2_{\geq 0}$ variable, this means that all three points~$(0, 0)$, $(K, 0)$, and~$(0, K)$ correspond to feasible assignments, but not~$(1, 1)$, see Figure~\ref{fig:xi-not-convex} for an illustration.
  Because~$(1, 1)$ is in the interior of the triangle mentioned above, we cannot separate it using only the~$\xi$ variable.
\end{example}
\begin{figure}[tb]
  \center
  \begin{tikzpicture}[scale=0.9]
    \draw[line width=3pt] (0, 2) -- (0, 0) -- (2,0) -- (2, 2);
    \draw (3, 0) rectangle node {$x_1$} (5, 1);
    \draw (3, 1) rectangle node {$x_2$} (5, 2);
    \draw (6, 0) rectangle node {$y_1$} (7, 2);
    \draw (7, 0) rectangle node {$y_2$} (8, 2);
    \draw (1,2.3) node {Bin};
    \draw (5.5, 2.3) node {Items};
  \end{tikzpicture}
  \hspace{1.5cm}
  \begin{tikzpicture}[scale=0.9]
    \draw[black!50] (-0.5, -0.5) grid (2.5, 2.5);
    \draw[-stealth, thick] (-0.5,0) -- (2.5, 0) node[right] {$\xi_x$};
    \draw[-stealth, thick] (0,-0.5) -- (0, 2.5) node[right] {$\xi_y$};
    \filldraw[col3] (0, 0) circle[radius=3pt];
    \filldraw[col3] (0, 1) circle[radius=3pt];
    \filldraw[col3] (0, 2) circle[radius=3pt];
    \filldraw[col3] (1, 0) circle[radius=3pt];
    \filldraw[col3] (2, 0) circle[radius=3pt];
    \draw[col1] (1, 1) circle[radius=3pt];
  \end{tikzpicture}
  \caption{
    Illustration of Example~\ref{example:xi-not-convex}.
    The partial assignment instance (left) has an infeasible solution from which the~$\xi$ variable (right) cannot be separated.
  }
    \label{fig:xi-not-convex}
\end{figure}

\paragraph{Application to Conflict Cuts}
We illustrate both techniques for the MKP.
Let~$\mathcal{G}$ be a partition of the set of items~$N$ such that, for every~$G \in \mathcal{G}$ and~$j_1,j_2 \in G$, we have both~$a_{j_1} = a_{j_2}$ as well as~$c_{i,j_1} = c_{i,j_2}$ for all~$i \in M$, i.e., items~$j_1$ and~$j_2$ have the same weight and objective value.
Then, a symmetry group of the MKP is~$\Pi = \bigotimes_{G\in \mathcal{G}} \mathcal{S}_{G}$, where~$\pi \in \Pi$ acts on a solution~$x$ of~\eqref{eq:naive-MK} by permuting the indices of items, i.e., $\pi(x)_{i,j} = x_{i,\pi^{-1}(j)}$.
Then, Proposition~\ref{prop:separate-benders} shows that symmetric Benders cuts can be separated efficiently.

Let~$(C,i) \in \mathcal{I}$ be an index of a cut from~\eqref{eq:multknapsack-BD}.
Then, the family of symmetric cuts is the set~$\{(\pi^{-1}(C),i) : \pi\in\Pi\}$, and it is uniquely determined by the number of elements selected per group~$G$.
We can thus characterize the entire family by a representative vector~$R(C) \in \Z_+^{\mathcal{G}}$, where~$R_G(C) = |C\cap G|$ for all~$G \in \mathcal{G}$.
The cuts of the family derived from~$(C,i)$ are then given by~$\sum_{G\in \mathcal{G}} \sum_{j\in G\cap C'} x_{ij} \leq |C'| - 1$ for all~$C' \subseteq N$ with~$R(C') = R(C)$.

In practice, to find a violated cut for a solution~$\fixed{x}$, one can iterate through all~$i \in M$, compute~$N_i = \{j \in N : \fixed{x}_{ij} = 1\}$, and check if there is~$(C,i)$ in the list of known cuts that has~$R(N_i) = R(C)$.
If such a pair exists, we can immediately derive a symmetric Benders cut separating~$\fixed{x}$ without solving the subproblem.
This idea can be enhanced by observing that, if~$(C,i) \in \mathcal{I}$, also~$(C',i) \in \mathcal{I}$ for all~$C'$ with~$C \subseteq C' \subseteq N$.
Hence, instead of checking~$R(N_i) = R(C)$, we can check whether there is~$(C,i) \in \mathcal{I}$ with~$R(N_i) \geq R(C)$ to find a violated symmetric cut.

\paragraph{Symmetry aware Lifting}
Recall that when computing the lifting coefficients for a given cut,
\begin{equation}
    \lift_{ij} \define |C| -1 - \max\brackets{
    \sum_{k\in C} x_{ik} + \sum_{k\in S} \lift_{ik} x_{ik} : x\in P_i,\ x_{ij} = 1
  }
  \tag{\ref{eq:sequence-lifting}}
\end{equation}
is dependent on the previously computed values of the entries in~$S$.
This implies the existence of up to~$\Oh{|N|!}$ many different possible lifted cuts.
If the symmetry group of the problem is of the form~$\Pi = \bigotimes_{G\in\mathcal{G}}\sym{G}$, one would hope to reduce the number of different liftings to at most~$\Oh{|\mathcal{G}|!}$.
In general, however, this intuitive idea is wrong because lifting coefficients~$\lift_{ij}$, $j \in G \in \mathcal{G}$ do not necessarily take the same value.
In fact, the difference between two lifting coefficients can be arbitrarily large as we show next.
\begin{example}
  \label{example:arbitrary-lifting-diff}
  Expanding on Example~\ref{example:xi-not-convex}, let~$x_1, \dots, x_{n_1}\in\brackets{0, 1}$ denote the indicator variables for the horizontal rectangles of dimension~$(K\times 1)$, and~$y_1, \dots, y_{n_2}\in\brackets{0, 1}$ for the vertical rectangles of dimension~$(1\times K)$.
  Any~$C\subseteq [n_1]$ with~$|C| = K + 1$ induces the cut~$\sum_{j\in C} x_j \leq |C| - 1$.
  Computing the lifting coefficient for~$y_1$ via Equation~\eqref{eq:sequence-lifting} results in~$\lift^y_1 = |C| - 1$, because forcing a vertical rectangle into the bin by setting~$y_1 = 1$ does not allow to add any horizontal rectangle.
  The partially lifted inequality is therefore~$\sum_{k \in C} x_{ik} + (|C| - 1)y_1 \leq |C| - 1$.
  Computing then the lifting coefficient for~$y_2$ results in~$\lift^y_2 = 0$, because fixing~$y_2 = 1$ allows to make use of the vertical rectangle associated with~$y_1$, achieving an objective value of~$|C| - 1$ in the maximization problem in~\eqref{eq:sequence-lifting}.
\end{example}
Since computing the lifting coefficients via~\eqref{eq:sequence-lifting} can be expensive, we suggest an alternative approach for computing weaker lifting coefficients~$\lift'_{ij} \leq \lift_{ij}$.
There is thus a trade-off between the strength of the cut and the computational cost.
\begin{definition}
  \label{def:replacement-coef}
  We consider a problem of the form~\eqref{eq:naive-MK}.
  Let~$(C, i)\in\mathcal{I}$ and~$j\notin C$.
  For any partial assignment~$\fixed{x}_i\in P_i$, with~$\fixed{x}_{ij}=1$, the \emph{replacement coefficient}~$\lift'_{ij}\in\Z_{\geq 0}$ is the largest integer such that any~$S\subseteq \brackets{k\in C: \fixed{x}_{ik} = 0}$ with~$|S|\leq \lift'_{ij}$ induces a feasible partial assignment~$\tilde{x}_i$ where~$\tilde{x}_{ij} = 0$, $\tilde{x}_{ik} = 1$ for all~$k\in S$ and~$\tilde{x}_{ik} = \fixed{x}_{ik}$ for all other entries.
\end{definition}
In other word, item~$j$ can always be ``replaced'' by any~$\lift'_{ij}$ elements from~$C$, which implies the following lemma.
\begin{lemma}
  Let~$(C, i)\in\mathcal{I}$ be a valid cut for~\eqref{eq:naive-MK}.
  Then, $\sum_{j\in C} x_{ij} + \sum_{j\notin C} \lift_{ij}' x_{ij} \leq \card{C}-1$ is a valid cut for~$P_i$.
\end{lemma}
\begin{proof}
  Let~$\fixed{x}_i\in P_i$.
  Since all~$\lift'$ are integer, let~$K\in\Z_{\geq 0}$ such that~$\sum_{j\in C} \fixed{x}_{ij} + \sum_{j\notin C} \lift'_{ij} \fixed{x}_{ij} =K$.
  In the following, we prove that $K \leq |C| - 1$, i.e., the inequality holds for all feasible solutions as claimed.
  We distinguish two cases.
  Let us first assume that~$\fixed{x}_{ij} = 0$ for all~$j\notin C$.
  Then the inequality reduces to~$\sum_{j\in C} \hat{x}_{ij} \leq |C|-1$, which trivially holds as~$\sum_{j\in C} x_{ij} \leq |C|-1$ is valid by assumption.
  It follows immediately that
  \begin{equation}
    K = \sum_{j\in C} \fixed{x}_{ij} + \sum_{j\notin C} \lift_j' \fixed{x}_{ij} = \sum_{j\in C} \fixed{x}_{ij} + 0 \leq |C| -1.
  \end{equation}
  Let us now assume that there exists a~$k\in [n]\setminus C$ with~$\fixed{x}_{ik} = 1$.
  Because~$\fixed{x}_i$ is a feasible solution, $C_0(\fixed{x}_i) \define\brackets{j\in C: \fixed{x}_{ij} = 0}$ is not empty and by definition of~$\lift'_{ik}$, $|C_0(\fixed{x}_i)|\geq \lift'_{ik}$.
  We can then create an equivalent feasible solution~$\tilde{x}_i\in P_i$ by replacing element~$k$ with~$\lift'_{ik}$ elements of~$C_0(\fixed{x}_i)$.
  This exchange preserves~$\sum_{j\in C} \tilde{x}_{ij} + \sum_{j\notin C} \lift'_{ij} \tilde{x}_{ij} =K$ but~$|C_0(\tilde{x}_i)| < |C_0(\fixed{x}_i)|$.
  One can then iterate this process until either~$C_0$ or~$\brackets{j\notin C : \tilde{x}_{ij} = 1}$ becomes empty.
  Because this process preserves feasibility, it cannot end with~$C_0 = \emptyset$.
  It then terminates with~$K = \sum_{j\in C} \tilde{x}_{ij} \leq |C|-1$ which proves the assertion.
\end{proof}
The advantages of this method are that the lifting coefficient is now sequence independent and its value is unique per item equivalence class.
Indeed, the construction does not take into account any previous lifting coefficient and since feasibility is preserved under symmetry, it follows that two symmetric variables replace equivalent sets.
In addition, one can weaken the construction a little more with the same definition, but limiting~$\lift^k_j = \min\{\lift'_j, k\}$.
Depending on the instance, $\lift^k_j$ can be found analytically for the lower values of~$k$.
The disadvantage, however, is that the cut becomes weaker as~$\lift^k_{ij} \leq \lift'_{ij} \leq \lift_{ij}$ for any~$k\in [n]$ and any set of previously lifted entries~$S\subset [n]\setminus C$ in~\eqref{eq:sequence-lifting}.

\begin{example}
  Continuing Example~\ref{example:multipleknapsack}, computing the lifting coefficient~$\lift_{ij}$ for a given knapsack~$i$ and item~$j\notin S\cup C$, one needed to solve a problem of the form
  \begin{equation}
    |C| - 1 - 
    \max\brackets{
      \sum_{k\in C}x_{ik} + \sum_{k\in S\setminus C} \lift_{ik} x_{ik} :
      \sum_{k\neq j} a_k x_{ik} \leq b_i - a_j
    },
  \end{equation}
  which is a (single) knapsack problem.
  Our lifting coefficients~$\lift'_{ij}$, however, are in this context much easier to compute. 
  Observe that $\lift'_{ij} \geq 1$ is equivalent to~$a_j \geq a_k$ for all~$k\in C$ and in general~$\lift'_{ij} = k$ is equivalent to~$a_j \geq a_{i_1} +\dots + a_{i_k}$ for any pairwise different~$i_1, \ldots, i_k \in C$.
  We can therefore reorder~$C = \brackets{i_1, \dots, i_{|C|}}$ such that~$a_{i_k} \geq a_{i_{k+1}}$ for all~$k\leq |C| -1$.
  The lifting coefficient is then simply~$\lift_{ij}' = \max\brackets{k : a_j \geq a_{i_1} + \dots + a_{i_k}}$.
  In the context of the knapsack problem, this is a well known lifting coefficient and it has been shown to be close to yielding facet-defining inequalities~\cite{balas1978facets,nemhauser1994lifted} (denoted~$\pi_j$ and~$\beta_j$ in the respective references).
\end{example}

\paragraph{Extended Formulation}
By introducing~$\zeta$-variables for each slot~$i \in M$, Inequality~\eqref{eq:generic-EF-cut} for a conflict cut~$(C,i)$ can be easily phrased as~$\sum_{G\in \mathcal{G}} \sum_{k = 1}^{|G\cap C|} \zeta^k_{iG} \leq |C| -1$.
Since this inequality is violated if and only if~$\zeta^{|G\cap C|}_{iG} = 1$ for all~$G \in \mathcal{G}$ (recall~$\zeta^k_{iG} = 1$ if and only if~$\sum_{j \in G}x_{ij} \geq k$), the inequality can be simplified to~$\sum_{G\in \mathcal{G}} \zeta_{iG}^{|C\cap G|} \leq |\mathcal{G}| -1$.
Alternatively, we can lift these inequalities with the replacement coefficients~$\lift'$ (or~$\lift^k$) as~$\lift'_{ij}$ is invariant for any~$j\in G\setminus C$.
The following proposition shows how these two can be applied in conjunction.
\begin{proposition}
  Let~$(C, i)\in\mathcal{I}$ be a valid cut for~\eqref{eq:naive-MK}.
  We assume that~$P_i$ has no trivially infeasible items, i.e., for every~$j\in N$, the partial assignment~$\tilde{x}_i$ where~$\tilde{x}_{ij} = 1$ and~$\tilde{x}_{ik} = 0$ for all other entries is feasible.  
  Let the symmetry group~$\Pi$ adhere to the conditions in Proposition~\ref{prop:separate-benders}.
  Denoting~$\lift'_{iG}$ the lifting coefficient of any element~$j\in G\setminus C$, the following are equivalent:
  \begin{align}
    \sum_{j\in C} \pi(x)_{ij} + 
    \sum_{j\notin C} \lift'_{ij} \pi(x)_{ij} &\leq |C|-1,
    & \forall\pi\in\Pi, \label{eq:lifted-vanilla}
    \\
    \sum_{G\in\mathcal{G}}\left(
      \sum_{k=1}^{|G\cap C|} \zeta^k_{iG} 
      +\sum_{k = |G\cap C| + 1}^{|G|} \lift'_{iG} \zeta^k_{iG}
    \right) &\leq |C| - 1.
    \label{eq:lifted-zeta}
  \end{align}
\end{proposition}
\begin{proof}
  Let the sequences~$\alpha^1_G, \dots \alpha^{|G|}_G$ for all~$G\in\mathcal{G}$ be defined as
  \begin{equation*}
    \alpha^k_G = \begin{cases}
      1, &\text{if } k\leq |G\cap C|, \\
      \lift'_{iG}, &\text{if } k > |G\cap C|.
    \end{cases}
  \end{equation*}
  We can then rewrite Equation~\eqref{eq:lifted-zeta} as~\eqref{eq:generic-EF-cut} and apply Theorem~\ref{thm:EF-works} if~$\alpha^k_G \geq \alpha^{k+1}_G$ for all~$k\leq |G|-1$ and all~$G\in \mathcal{G}$.
  Observe in particular that if either~$G\cap C$ or~$G\setminus C = \emptyset$, $\alpha_G$ is constant and as such trivially sorted.
  It suffices then to prove that~$\lift'_{iG} \leq 1$ if~$C\cap G \neq \emptyset$.
  
  Let the pair~$j, j'\in G$ be such that~$j\notin C$ and~$j'\in C$ and assume by contradiction that~$\lift_{iG}' \geq 2$.
  Let~$\fixed{x}_i$ be the feasible partial assignment where~$\fixed{x}_{ij} = 1$ and~$\fixed{x}_{ik} = 0$ for any~$k\in N\setminus\brackets{j}$.
  By the definition of~$\lift'_G$, we can create a feasible solution by swapping~$\fixed{x}_{ij}$ with any two elements from~$C_0(\fixed{x}_i)\define\brackets{k\in C : \fixed{x}_{ik}= 0}$.
  In particular, let one of the two elements be~$j'$.
  We then get the feasible solution~$x'_i\in P_i$ with~$x'_{ij} = 0$, $x'_{ij'} = 1$ but~$\left|C_0(x'_i)\right| \leq \left|C_0(\fixed{x}_i)\right|-2$. 
  Because of~$\Pi$'s structure, there exists a transposition~$\gamma\in\Pi$ that swaps~$j$ and~$j'$ since they both belong to the same~$G$.
  We then have that~$x^*_i\define \gamma(x'_i)\in P_i$ is a feasible partial solution with~$x^*_{ij} =1$, $x^*_{ij'} = 0$ but it uses one more element in~$C$ than~$\fixed{x}_i$ as~$|C_0(x^*_i)| \leq |C_0(\fixed{x}_i)| - 1$.
  Repeating this argument, we can iteratively reduce~$|C_0|$ by 1 until~$|C_0| =1$.
  Denoting this assignment by~$x^\star_i$, we have in particular that~$C_0(x^\star_i) = \brackets{j'}$ and~$x^\star_{ij} = 1$.
  By symmetry, $\gamma(x^\star_i)$ must then also be a feasible partial assignment, however~$C_0(\gamma(x^\star_i)) = \emptyset$ which is a contradiction because no feasible solution can assign all elements in~$C$ to slot~$i$.
\end{proof}
\begin{remark}
  Note that in the degenerate case where there is some equivalence class~$G$ with~$x_{i}\notin P_i$ if~$\sum_{j\in G} x_{ij} > 0$, the lifting coefficient would be~$\lift'_{ij} = |C|$ for all~$j\in G$, independently if there are equivalent items in~$C$ or not.
  It follows then that for this group the sequence~$\alpha_G$ should instead be~$\alpha^k_G = \lift'_{iG}$ if~$k\leq |G\setminus C|$ and~$\alpha^k_G = 1$ afterwards, and the corresponding sum in~\eqref{eq:lifted-zeta} has then the form~$\sum_{k=1}^{|G\setminus C|} \lift'_{iG} \zeta^k_{iG} +\sum_{k = |G\setminus C| + 1}^{|G|} \zeta^k_{iG}$.

\end{remark}

%
%

\section{Applications}\label{sec:Applications}
The discussion of the previous section focused on general techniques to exploit symmetry in BD.
In the following, we illustrate how these techniques can be applied for three different applications with varying difficulty of the cut generation subproblem:
sphere packing (\spherepack), rectangle packing (\rectanglepack), and machine scheduling (\scheduling).
In the following, we discuss how the ideas explained in Secs.~\ref{sec:Benders-Symmetry} and~\ref{sec:usingsymmetry} can be realized for these examples.

%
%

\medskip
\noindent\textbf{2D Bin Packing}
Given a set of items~$N$ and bins~$B$, the bin packing problem is to find the least number of bins that can hold all the items.
We assume that all bins~$i \in B$ are rectangles of width~$W_i$ and height~$H_i$, and consider two different types of items: spheres and rectangles.
For both variants, we use the same master problem~\eqref{eq:BP-master-problem}.
Note that, as this is an assignment problem, its formulation stems from the generic one in~\eqref{eq:naive-MK}.
\begin{subequations}
    \label{eq:BP-master-problem}
    \begin{align}
        \min \sum_{i\in B} u_i & &&\nonumber \\
        \sum_{i\in B} x_{ij} &= 1, && j\in N, \label{eq:master-coveringcons}\\
        \sum_{j\in C} x_{ij} &\leq |C| - 1, &&  (C, i)\in \mathcal{I}, \label{eq:master-benders-cut} \\
        \sum_{j \in N} a_j x_{ij} &\leq W_i H_i u_{i}, && i \in B, \label{eq:master-areacons} \\
        x_{ij} &\leq u_i, && j\in N,\;  i \in B. \label{eq:master-openbincons}
    \end{align}    
\end{subequations}
In addition to the constraints from~\eqref{eq:naive-MK}, the master problem adds auxiliary variables~$u_i\in \brackets{0, 1}$, $i \in B$, indicating whether the corresponding bin~$i$ has assigned items.
The coefficients~$a_j$, $j \in N$, correspond to the area of item~$j$.
Thus, Equation~\eqref{eq:master-areacons} strengthens the master problem~\eqref{eq:naive-MK} by excluding all solutions that assign a single bin~$i \in B$ items whose total area exceeds the bin's area~$W_iH_i$.

Given a solution~$\fixed{x}$, the Benders subproblem is to decide, for each bin~$i \in B$, if the selected items~$\{j \in N: \fixed{x}_{ij} = 1\}$ fit into the bin.
For \spherepack, this problem can be solved by deciding feasibility of the following system,  cf.~\cite{khajavirad2024circle}:
\begin{subequations}
  \label{eq:subproblem-2DSP}
  \begin{align}
    \left(y^{(1)}_j - y^{(1)}_k\right)^2 + \left(y^{(2)}_j - y^{(2)}_k\right)^2 &\geq (\fixed{x}_{ij} + \fixed{x}_{ik} - 1)(r_j + r_k)^2, && j, k\in N,\; j \neq k, \label{eq:2DSP-non-overlap}\\
    r_j \leq y^{(1)}_j &\leq W_i - r_j, && j\in N, \\
    r_j \leq y^{(2)}_j &\leq H_i - r_j, && j\in N.
\end{align}
\end{subequations}
The variable pair~$\left(y^{(1)}_j, y^{(2)}_j\right)$ acts as the position of the center of sphere~$j$ in the bin.
Observe that Equation~\eqref{eq:2DSP-non-overlap} enforces that two spheres do not overlap only if they are both assigned to bin~$i$ ($\fixed{x}_{ij} = \fixed{x}_{ik} = 1$).
To build an SDG for Subproblem~\eqref{eq:subproblem-2DSP}, one can either use rules for general MINLPs, see~\cite{liberti2012reformulations}, or capture the symmetries by a more compact tailored graph.
To this end, recall that the SDG only needs to capture the symmetries of the~$x$-variables, because the~$y^{(1)}$ and~$y^{(2)}$ variables are not present in the master problem.
More precisely, we define the node set
\[
  V = \{a_i\} \cup \bigcup_{j \in N} \{v_{ij}, r_j\},
\]
where~$a$ is the anchor of the graph, $v_{ij}$ is a variable node representing variable~$x_{ij}$, and~$r_j$ represents the radius of sphere~$j$.
The anchor~$a_i$ receives a color that uniquely characterizes the dimensions~$(W_i,H_i)$ of the respective bin~$i$, the variable nodes are colored according to their objective coefficient and type in the master problem, and the node~$r_j$ receives a unique color that indicates it as radius-node.
Moreover, we define the edge set
\[
  E = \bigcup_{j \in N} \brackets{\{a_i,r_j\}, \{r_j,v_{ij}\}},
\]
where~$\{a_i,r_j\}$ receives a color corresponding to the radius~$r_j$, and edge~$\{r_j,v_{ij}\}$ remains uncolored.
Due to this construction, an automorphism of~$G$ can only exchange variable nodes~$v_{ij}$ and~$v_{ij'}$ if~$r_j = r_{j'}$.
Thus, the SDG~$G$ captures the symmetries of the subproblem.

If the subproblem~\eqref{eq:subproblem-2DSP} is determined to be infeasible, the corresponding cut is added to~\eqref{eq:master-benders-cut}.
Said cut can be strengthen via the lifting technique described earlier:
In the context of \spherepack, determining the value of~$\lift'_{ij}$ for a given pair~$(C, i)\in \mathcal{I}$ and~$j\in N\setminus C$ corresponds finding the largest~$k$ such that any subset of spheres~$S\subseteq C$ of cardinality~$|S| = k$ can fit inside a sphere of radius~$r_j$.
Since this problem can be challenging, we restrict ourself to finding~$\lift^2 = \min\brackets{\lift', 2}$.
In this case, $\lift^2_j \geq 1$ is equivalent to the sphere~$j$ being larger than any sphere in~$C$, and~$\lift^2_j = 2$ means that it can even cover any two non-overlapping spheres from~$C$.
This value can be analytically derived by simply reordering~$C = \brackets{j_1, ..., j_{|C|}}$ such that~$r_{j_k} \geq r_{j_{k+1}}$ for all~$k\leq |C| -1$.
It is easy to see that
\begin{equation}
  \lift^2_{ij} = \begin{cases}
    2, & \text{if } r_j \geq r_{j_1} + r_{j_2}, \\
    1, & \text{if } r_j \in [r_{j_1}, r_{j_1} + r_{j_2} ), \\
    0, & \text{if } r_j < r_{j_1}.
  \end{cases}
\end{equation}

The subproblem of~\rectanglepack can be modeled as a MIP, cf.~\cite{pisinger2007using}, and an SDG can either be derived from this MIP, or via a similar construction as for~\spherepack.
See Appendix~\ref{app:rectangle-packing} for further details.

%
%

\medskip
\textbf{Machine Scheduling With Setup Times}
We consider a machine scheduling problem with setup times between different jobs.
We are given a set of machines~$M$, set of jobs~$N$, processing time~$p_{ij}$ for job~$j$ on machine~$i$, and setup times~$s_{ijk}$ for processing job~$k$ directly after job~$j$ on machine~$i$ that satisfy the triangle inequality~$s_{ijk} + s_{ikl} \geq s_{ijl}$.
The objective is to find the minimum required time to have all jobs processed, called the \emph{makespan}, given that each machine can only process one job at a time.
Using the LBBD approach of~\cite{tran2016decomposition}, we add to~\eqref{eq:naive-MK} a variable~$T \geq 0$ representing the makespan, and set the objective to minimize it.
To link~$T$ with the remaining variables, \cite{tran2016decomposition} introduces Benders cuts
\begin{equation}
  T \geq T_i(C) - \sum_{j\in C} (1-x_{ij})\theta_{ij}, \quad (C, i)\in \mathcal{I},
  \label{eq:scheduling-BD}
\end{equation}
where~$C \subseteq N$, $i \in M$, and~$T_i(C)$ denotes the minimum makespan for processing all jobs in~$C$ on machine~$i$.
The coefficients~$\theta_{ij}\define p_{ij} + \max\brackets{s_{ikj} : k\in C}$ represent an upper bound on the time saved when not assigning job~$j$ to that machine.
Next to~\eqref{eq:scheduling-BD}, \cite{tran2016decomposition} strengthen the model by introducing auxiliary continuous variables.
Since the details are not relevant for the following discussion, we refer to Appendix~\ref{app:MS} for the finer technicalities.

\begin{figure}[t]
  \centering
  \begin{tikzpicture}[scale=0.75]
    \tikzstyle{n} += [circle, draw=black];
    \node (z_1) at (  1,   1) [n, fill=col3] {};
    \node (z_2) at (3.5,   1) [n, fill=col4] {};
    \node (z_3) at (3.5, 3.5) [n, fill=col4] {};
    \node (z_4) at (  1, 3.5) [n, fill=col3] {};

    \draw (z_1) node[below left =2pt] {$x_1$};
    \draw (z_2) node[below right=2pt] {$x_2$};
    \draw (z_3) node[above right=2pt] {$x_3$};
    \draw (z_4) node[above left =2pt] {$x_4$};

    \begin{scope}[-stealth, thick]
      \draw[col6] (z_1) to[out=10,  in=170] (z_2);
      \draw[col6] (z_2) to[out=190, in=-10] (z_1);

      \draw[col1] (z_1) to[out=55,  in=215] (z_3);
      \draw[col1] (z_3) to[out=235, in=35 ] (z_1);

      \draw[col1] (z_1) to[out=100, in=260] (z_4);
      \draw[col5] (z_4) to[out=280, in=80 ] (z_1);

      \draw[col6] (z_2) to[out=100, in=260] (z_3);
      \draw[col5] (z_3) to[out=-80, in=80 ] (z_2);

      \draw[col1] (z_2) to[out=125, in=325] (z_4);
      \draw[col1] (z_4) to[out=-55, in=145] (z_2);

      \draw[col5] (z_3) to[out=170, in=10 ] (z_4);
      \draw[col5] (z_4) to[out=-10, in=190] (z_3);
    \end{scope}
  \end{tikzpicture}
  \hspace{1cm}
  \begin{tikzpicture}[scale=0.75]
    \tikzstyle{v} += [thick, -stealth];
    \tikzstyle{nb} += [circle, draw=black];
    \tikzstyle{nw} += [circle, draw=white];
    \node (x_1)    at (  2, 1) [nb, fill=col3] {};
    \node (x_2)    at (  3, 1) [nb, fill=col4] {};
    \node (x_3)    at (  4, 1) [nb, fill=col4] {};
    \node (x_4)    at (  5, 1) [nb, fill=col3] {};
    \node (e_{12}) at (  1, 3) [nw, fill=black] {};
    \node (e_{13}) at (  2, 3) [nw, fill=black] {};
    \node (e_{14}) at (  3, 3) [nw, fill=black] {};
    \node (e_{23}) at (  4, 3) [nw, fill=black] {};
    \node (e_{24}) at (  5, 3) [nw, fill=black] {};
    \node (e_{34}) at (  6, 3) [nw, fill=black] {};
    \node (c)      at (3.5, 4) [nw, fill=black!50] {};

    \begin{scope}[thick]
      \draw[col6] (x_1) -- (e_{12});
      \draw[col6] (x_2) -- (e_{12});
      \draw[col1] (x_1) -- (e_{13});
      \draw[col1] (x_3) -- (e_{13});
      \draw[col1] (x_1) -- (e_{14});
      \draw[col1] (x_2) -- (e_{24});
      \draw[col1] (x_4) -- (e_{24});
      \draw[col5] (x_4) -- (e_{14});
      \draw[col6] (x_2) -- (e_{23});
      \draw[col5] (x_3) -- (e_{23});
      \draw[col5] (x_3) -- (e_{34});
      \draw[col5] (x_4) -- (e_{34});
      \draw[black!40] (e_{12}) -- (c);
      \draw[black!40] (e_{13}) -- (c);
      \draw[black!40] (e_{14}) -- (c);
      \draw[black!40] (e_{23}) -- (c);
      \draw[black!40] (e_{24}) -- (c);
      \draw[black!40] (e_{34}) -- (c);      
    \end{scope}

    \foreach \n in {x_1, x_2, x_3, x_4}
        \draw[black] (\n) node[below=3pt]  {$\n$};

  \end{tikzpicture}
  \caption{
    Example of compact SDG for TSP. 
    Left: exemplary instance of the \scheduling subproblem with four selected jobs;
    node colors represent processing times, and arc colors the setup times.
    Right: corresponding SDG.
    One only needs to add a dummy node for each job pair and the grey anchor node, where~$x_1,\dots,x_4$ are variable nodes.
  }
  \label{fig:TSP-SDG}
\end{figure}

Given a solution~$\fixed{x} \in \B{M \times N}$, the subproblem is to determine, for all~${i \in M}$, the value~$T_i(N_i)$, where~$N_i = \{j \in N : \fixed{x}_{ij} = 1\}$.
Note that for any subset of jobs~$C$, the minimum makespan is given by~$T_i(C) = \sum_{j\in C} p_{ij} + S_i(C)$, where~$S_i(C)$ is the minimum total setup time of jobs in~$C$ on machine~$i$.
Since the setup times depend on the ordering, finding~$S_i(C)$ can be reduced to solving a traveling salesperson problem (TSP), see~\cite{tran2016decomposition}.
Note that the classical subtour elimination formulation of the TSP has exponentially many constraints.
To detect symmetries, standard approaches list all these constraints explicitly, which is impractical.
Theorem~\ref{thm:merge-SDGs}, however, allows to define a compact SDG for the subproblem.
Let the nodes of the machine scheduling problem be given by
\begin{equation*}
  V = \brackets{a} \cup \bigcup_{j\in N} v_j \cup \bigcup_{\brackets{j, k} \subseteq N} u_{jk}.
\end{equation*}
In particular, we can create variable nodes~$v_j$ for each job (respectively city, in the context of TSP) $j$ and uncolored auxiliary nodes~$u_{jk}$ for any pair~$\brackets{j, k}\subseteq N$.
We color the variable nodes according to their processing time and connect node~$v_j$ to~$u_{jk}$ with an edge whose color correspond to the setup time (travel time) from~$j$ to~$k$.
See Fig.~\ref{fig:TSP-SDG} for an illustration.

\paragraph{Cut Pool and Extended Formulation}
Since solving the TSP subproblem can be costly, the symmetric cut pool can potentially save a lot of time.
When scanning the cut pool, we need to decide if the set of jobs~$C$ used in a previously derived cut~$(C,i)$ is symmetric to the current assignment~$N_i = \{j \in N : \fixed{x}_{ij} = 1\}$, $i \in M$.
To this end, we define for each machine~$i \in M$ an equivalence relation~$\sim_i$ defining the groups of identical jobs.
We say that two jobs~$j_1, j_2 \in N$ satisfy
\begin{equation}
  j_1\sim_{i} j_2 \iff 
  \begin{cases}
    &p_{i j_1\phantom{k}} = p_{i j_2},\ s_{i j_1 j_2} = s_{i j_2 j_1}, \\
    &s_{i j_1 k} = s_{i j_2 k}, \text{ for all } k\in N\setminus\{j_1, j_2\}, \\
    &s_{i k j_1} = s_{i k j_2}, \text{ for all } k\in N\setminus\{j_1, j_2\}.
  \end{cases} \nonumber
\end{equation}
This relation induces job groups~$G\in\mathcal{G}_i$ that depend on machine~$i \in M$, and again, we can define a representative vector~$R(N_i) \in \Z_+^{\mathcal{G}_i}$, $R_G(N_i) = |N_i \cap G|$, to decide whether a symmetric violated cut exists in the pool.
Moreover, by introducing~$\zeta$-variables, all cuts symmetric to~\eqref{eq:scheduling-BD}, can be represented by
\begin{equation}
  T \geq T_i(C) - \sum_{G\in \mathcal{G}_i} \sum_{h = 1}^{|G\cap C|} (1 - \zeta^h_{iG})\theta_{iG}, \quad (C, i)\in \mathcal{I}\label{eq:EF-scheduling-cut}
\end{equation}
where~$\theta_{iG} \define p_{iG} + \max\{s_{i H G} : H \in \mathcal{G}_i,\ H \cap C \neq \emptyset\}$,
$p_{iG} = p_{ij}$ and $s_{iGH} = s_{ijk}$ for all~$j\in G$, $k\in H$ and groups~$G,H \in \mathcal{G}_i$.

%
%

\section{Numerical Results}
\label{sec:numerics}

Sec.~\ref{sec:usingsymmetry} has introduced different means to exploit symmetry in BD.
This section's goal is to empirically evaluate the impact of these methods.
We are particularly interested in the impact of state-of-the-art symmetry handling methods when applied within BD, and if the cut pool or extended formulation has a stronger impact on the performance of BD.
The former has not been investigated yet because symmetry detection for BD was not available;
regarding the latter, although cut pools and extended formulations have been discussed before~\cite{Daryalal2023,Han2025,riise2016recursive}, to the best of our knowledge, they have not been compared with each other.

\paragraph{Implementation Details}
We have implemented LBBD schemes for the \rectanglepack, \spherepack, and \scheduling problem using the solver \scip~\cite{SCIP8}.
For each problem, we implemented a so-called constraint handler to separate Benders cuts.
To enable automatic symmetry detection for BD in \scip, we implemented callbacks of the respective constraint handlers that add symmetry information of the Benders oracles to \scip's SDG.
Our source code is available at GitHub\footnote{\url{https://github.com/Cedric-Roy/Benders-Symmetry-Supplements}} and zenodo~\cite{releaseBenders}.

Symmetric cut pools are implemented by maintaining a list of abstract Benders cuts that have already been separated.
An \emph{abstract cut} does not make use of explicit variables (items/jobs in our applications), but only counts how many variables of which symmetry class~$G\in\mathcal{G}$ have been involved in the cut via representative vectors~$R(C)$.
Abstract cuts thus correspond to families of cuts.
Whenever a solution needs to be separated, we iterate through the abstract cuts in our list.
For each abstract cut, we check whether a member of its associated family of cuts is violated using Proposition~\ref{prop:separate-benders}.
Only when no violated cut is found in the list, the separation oracles are triggered to possibly generate a new cut, which is then stored in our list of cuts.

Extended formulations are implemented by adding the auxiliary~$\zeta$-variables to the initial Benders master problem.
The separation of cuts then follows the same framework as for the cut pool, but cuts are added in terms of the~$\zeta$-variables, as indicated by Theorem~\ref{thm:EF-works}.
For both the cut pool and extended formulation, we only separate integer solutions to avoid too many calls of the separation oracles.

\medskip
\noindent\emph{Settings}
We compare the performance of a default implementation of LBBD, not exploiting symmetry information (setting \setting{plain}) with LBBD making use of symmetric cut pools (\setting{pool}) and extended formulations.
For extended formulations, we investigate the effect of adding Benders cuts as constraints to the problem (\setting{EFcons}) or as rows to the LP relaxation (\setting{EFrow}).
Constraints remain in the problem, whereas rows can be removed again.
Using rows, the LP relaxation has potentially less constraints and a weaker bound, but can be solved potentially faster.
We also investigate the effect of enabling/disabling \scip's built-in default symmetry handling methods.
Additionally, we investigate the effect of the lifted conflict cuts for \spherepack and \rectanglepack.
As such, we created an additional batch of settings to evaluate our test set over, namely: \setting{lifted}, \setting{lifted pool}, and \setting{lifted EFcons} which correspond to lifted cuts of \setting{plain}, \setting{pool}, and \setting{NC EFcons}, respectively.
Note that there are two equivalent ways of denoting the cuts in the extended formulation: the non-compact (NC) and the compact version
\begin{equation}
  \label{eq:non-compact-vs-compact}
  \setting{NC EFcons}: \sum_{G\in\mathcal{G}} \sum_{k = 1}^{|G\cap C|} \zeta_{iG}^k \leq |C| - 1, \quad
  \setting{EFcons}: \sum_{G\in \mathcal{G}} \zeta_{iG}^{|C\cap G|} \leq |\mathcal{G}| -1.
\end{equation}
Only the non-compact formulation is compatible with the lifting procedure discussed previously.
This means that we can only evaluate the benefit from lifting with the difference between \setting{lifted EFcons} and \setting{NC EFcons}.

All experiments use a time limit of \SI{2}{\hour}.
Mean numbers are reported in shifted geometric mean~$\prod_{i=1}^n (t_i + s)^{\nicefrac{1}{n}} - s$ for numbers~$t_1,\dots,t_n$.
We use~$s=10$ for nodes of the branch-and-bound tree, and~$s=1$ for all remaining quantities.

\medskip
\noindent\emph{Hardware and Software Specifications}
\label{app:hardware}
All experiments have been conducted on a Linux Cluster with Intel Xeon E5-1620 v4 \SI{3.5}{\giga\hertz} quad core processors and \SI{32}{\giga\byte} memory.
The code was executed single-threaded with a time limit of \SI{2}{\hour}.
We use \scip~9.2.3~\cite{SCIP9} as branch-and-bound framework; LP relaxations are solved using \solver{SoPlex}~7.1.5 and nonlinear problems are solved by \solver{Ipopt}~3.14.20~\cite{WachterBiegler2006}.
Symmetries of SDGs are detected by \solver{sassy}~1.1~\cite{AndersEtAl2023} and \solver{Nauty}~2.8.8~\cite{Nauty}.

\paragraph{Test Sets}
We have generated random instances of the \rectanglepack, \spherepack, and \scheduling problem as follows.

\medskip
\noindent\emph{Sphere Packing}
All instances make use of five rectangular bins of width~$W = 15$ and height~$H = 10$.
The instances differ, however, in the number of spheres and their radii.
To guarantee that the instances have symmetries, we create spheres in batches of identical radii.
Each batch has size five, and for each batch we sample its radius uniformly at random from the interval~$[1,3]$ and round it to one decimal.
We have created our instances by either sampling two or three batches, i.e., our instances have either~$10$ or~$15$ spheres.
For each number of batches, we have created~$40$ random instances, thus yielding~$80$ instances in total.

\medskip
\noindent\emph{Rectangle Packing}
The instances of the rectangle packing problem are created similarly.
Instead of making use of five bins, however, we create one bin per item, and bins have width~$W = 15$ and height~$H = 15$.
Moreover, for each instance, we create five batches consisting of either four or five items.
The width and height of each batch of items is sampled uniformly at random from the interval~$[1,10]$ and is rounded to one decimal.
For each number of batches, we have created~$40$ random instances, thus yielding~$80$ instances in total.

\medskip
\noindent\emph{Machine Scheduling}
The instances of the machine scheduling problem have been generated similarly, except for the setup times.
Each instance consists of three identical machines.
We divided the instances in three batches each consisting of five, six, or seven jobs.
By adding the artificial job~$0$, this amounts in instances with~$16$, $19,$ or~$22$ jobs.
The processing time of each batch of jobs is sampled uniformly at random from the interval~$[1, 10]$ and is rounded to one decimal.

In order to generate setup times between job pairs that follow the triangle inequality, we took inspiration from the tool switching problem~\cite{tang1988models}.
For each job~$j$, we select a set~$C_j$ of tools, and the setup time between jobs~$j$ and~$j'$ is computed as the cardinality of~$(C_j \setminus C_{j'}) \cup (C_{j'} \setminus C_j)$, i.e., the symmetric difference of~$C_j$ and~$C_{j'}$.
This notion of setup time automatically satisfies the triangle inequality.
In our implementation, the sets of tools are randomly selected subsets~$C\subseteq\brackets{1, \dots, 10}$ with~$|C|\leq 5$.
Finally, we force the setup times~$s_{ij0} = 0$ for all jobs~$j$ and all machines~$i$. 
For each type listed above, we generated batches of~$40$ instances, thus yielding~$120$ instances in total.

\paragraph{Numerical Results}
The results are summarized in Table~\ref{tab:statistics}, where column ``\#solved'' reports on the number of solved instances, ``time'' and ``cut-time'' give the mean solving time and time needed to separate Benders cuts, respectively; ``\#sepa'' is the mean number of solutions that have been separated.

\begin{table}[t]
  \caption{Statistics on different approaches for solving \rectanglepack, \spherepack, and \scheduling problems.}
  \label{tab:statistics}
  \centering
  \begin{scriptsize}
  \begin{tabular*}{\textwidth}{@{}l@{\;\;\extracolsep{\fill}}rrrrrrrr@{}}
    \toprule
    & \multicolumn{4}{c}{with symmetry handling} & \multicolumn{4}{c}{without symmetry handling}\\
    \cmidrule{2-5} \cmidrule{6-9}
    setting & \#solved & time & cut-time & \#sepa & \#solved & time & cut-time & \#sepa \\
    \midrule
    \multicolumn{9}{@{}l}{\rectanglepack:}\\
     plain & \num{ 24} & \num{2678.0} & \num{2613.6} & \num{     894.4} & \num{ 24} & \num{2704.8} & \num{2642.1} & \num{    1022.0} \\
      pool & \num{ 35} & \num{1159.9} & \num{ 395.0} & \num{    2600.1} & \num{ 26} & \num{1960.4} & \num{ 498.5} & \num{    6981.4} \\
     EFrow & \num{ 42} & \num{1059.5} & \num{ 935.7} & \num{     979.1} & \num{ 34} & \num{1451.0} & \num{1286.6} & \num{    1480.8} \\
    EFcons & \num{ 59} & \num{ 232.9} & \num{ 118.4} & \num{     103.7} & \num{ 41} & \num{ 793.3} & \num{ 233.7} & \num{     199.0} \\
    \midrule
    \multicolumn{9}{@{}l}{\spherepack:}\\
     plain & \num{ 46} & \num{ 287.4} & \num{ 287.1} & \num{       8.8} & \num{ 42} & \num{ 357.9} & \num{ 357.8} & \num{      11.6} \\
      pool & \num{ 57} & \num{ 109.2} & \num{ 109.0} & \num{      10.1} & \num{ 55} & \num{ 127.8} & \num{ 113.3} & \num{      24.3} \\
     EFrow & \num{ 55} & \num{ 202.7} & \num{ 202.4} & \num{       7.1} & \num{ 46} & \num{ 250.0} & \num{ 249.8} & \num{       9.0} \\
    EFcons & \num{ 55} & \num{ 173.4} & \num{ 173.1} & \num{       5.7} & \num{ 53} & \num{ 187.4} & \num{ 187.0} & \num{       6.6} \\
    \midrule
    \multicolumn{9}{@{}l}{\scheduling:}\\
     plain & \num{  7} & \num{6801.9} & \num{6274.9} & \num{  230070.3} & \num{  8} & \num{6790.8} & \num{6263.8} & \num{  229783.5} \\
      pool & \num{ 45} & \num{3318.7} & \num{  29.7} & \num{ 2018389.5} & \num{ 46} & \num{3286.2} & \num{  30.2} & \num{ 2039326.9} \\
     EFrow & \num{ 84} & \num{ 561.7} & \num{   6.4} & \num{  143804.3} & \num{ 83} & \num{ 588.6} & \num{   6.8} & \num{  158353.7} \\
    EFcons & \num{120} & \num{  29.9} & \num{   1.0} & \num{     126.7} & \num{119} & \num{  30.5} & \num{   1.1} & \num{     145.5} \\
    \bottomrule
  \end{tabular*}
  \end{scriptsize}
\end{table}

In the introduction, we have posed two questions, which we now answer in turn.
Question~\eqref{Q1} concerned whether symmetries can be automatically detected when using BD.
Using the framework of SDGs, this is indeed possible and Table~\ref{tab:statistics} shows that we can greatly benefit from handling symmetries in BD when using \scip's default symmetry handling methods.
For the \rectanglepack problem, we can solve considerably more instances when handling symmetries and substantially reduce the running time (between~\SI{27}{\percent} and~\SI{71}{\percent} for \setting{pool}, \setting{EFrow}, and \setting{EFcons}).
Using \setting{plain}, however, the effect is moderate.
This indicates that symmetry handling only becomes effective when it is combined with exploiting symmetries in generating Benders cuts.

The effect of symmetry handling for the \spherepack problem is less pronounced with running time improvements between~\SI{8}{\percent} and~\SI{19}{\percent}.
To explain this behavior, we compared the size of the branch-and-bound trees for \setting{pool}, \setting{EFrow}, and \setting{EFcons} for the instances that could be solved by all three methods within the time limit.
The mean number of nodes in these trees is~11.5--16.5, see Appendix~\ref{app:detailedresults}, i.e., the effect of symmetry handling is minor due to the rather small trees.

Finally, symmetry handling seems to have almost no effect on the \scheduling problem, which can be explained by the symmetry handling methods used by \scip.
Since all master problem variables for the \rectanglepack and \spherepack problems are binary, \scip handles symmetries by the methods orbital fixing~\cite{OstrowskiEtAl2011,PfetschRehn2019} and lexicographic reduction~\cite{DoornmalenHojny2024a}.
The \scheduling problem, however, has predominantly continuous variables in the master problem, and \scip decides that it is more important to handle symmetries on the continuous variables using SST cuts~\cite{LibertiOstrowski2014,Salvagnin2018}.
SST cuts handle only a small amount of symmetries and their structure is comparable to inequalities that one would add by hand.
We therefore also tested whether the more sophisticated techniques lexicographic reduction and orbital reduction~\cite{DoornmalenHojny2024a} yield a performance improvement for \scheduling problems.
This can be easily achieved in \scip by setting the parameter \setting{misc/usesymmetry} to~$3$, see Table~\ref{tab:statisticssched} for details.
Indeed, using these techniques, \setting{plain} solves~61 more instances and reduces the running time by more than~\SI{92}{\percent}, and also \setting{pool}, \setting{EFrow}, and \setting{EFcons} substantially improve their performance.
In particular, the most competitive setting \setting{EFcons} reduces its running time by another~\SI{78}{\percent}.
We thus conclude that automatic symmetry detection and handling in BD is very important, because it allows solvers to make use of sophisticated (already built-in) symmetry handling techniques that users can not easily implement on their own.

\begin{table}[tb]
  \caption{Statistics on different approaches for solving \scheduling problems with default and advanced symmetry handling methods.}
  \label{tab:statisticssched}
  \centering
  \begin{scriptsize}
  \begin{tabular*}{\textwidth}{@{}l@{\;\;\extracolsep{\fill}}rrrrrrrr@{}}
    \toprule
    & \multicolumn{4}{c}{default symmetry handling} & \multicolumn{4}{c}{advanced symmetry handling}\\
    \cmidrule{2-5} \cmidrule{6-9}
    setting & \#solved & time & cut-time & \#sepa & \#solved & time & cut-time & \#sepa \\
    \midrule
     plain & \num{  7} & \num{6801.9} & \num{6274.9} & \num{  230070.3} & \num{ 68} & \num{ 511.6} & \num{ 455.8} & \num{   16226.3} \\
      pool & \num{ 45} & \num{3318.7} & \num{  29.7} & \num{ 2018389.5} & \num{ 76} & \num{ 160.6} & \num{   7.9} & \num{   56406.8} \\
     EFrow & \num{ 84} & \num{ 561.7} & \num{   6.4} & \num{  143804.3} & \num{ 99} & \num{  65.8} & \num{   3.6} & \num{   11828.9} \\
    EFcons & \num{120} & \num{  29.9} & \num{   1.0} & \num{     126.7} & \num{120} & \num{   6.5} & \num{   1.0} & \num{      86.1} \\
    \bottomrule
  \end{tabular*}
  \end{scriptsize}
\end{table}

Question~\eqref{Q2} asked whether the approach using a cut pool or extended formulation is better suited to exploit symmetries in BD.
In the following discussion, we focus on the results for BD that handles symmetries of the master problem.
From Table~\ref{tab:statistics} it is apparent that both approaches substantially improve \setting{plain} BD, reducing the mean running times for \rectanglepack, \spherepack, and \scheduling up to~\SI{91}{\percent}, \SI{62}{\percent}, and~\SI{99}{\percent}, respectively.
In particular, while \setting{plain} could only solve~7 of the~120 \scheduling instances, \setting{EFcons} solves all instances.

For the \rectanglepack and \scheduling problem, \setting{EFcons} clearly dominates \setting{pool} and \setting{EFrow} when separating solutions.
A possible explanation is that adding Benders cuts as constraints in an extended formulation guarantees that no symmetric solutions can be computed anymore.
The search space for \setting{EFcons} is thus considerably smaller than for \setting{EFrow} and \setting{pool}.
This hypothesis can be confirmed by comparing the number of nodes in the branch-and-bound trees of the instances that are solved by all three settings, see Appendix~\ref{app:detailedresults}:
For the three test sets, the mean trees for \setting{EFcons} are 30--\SI{95}{\percent} smaller than for \setting{EFrow} (resp.\ 22--\SI{99}{\percent} for \setting{pool}).

For \spherepack problems, however, \setting{EFcons} is considerably slower than \setting{pool}, which has two explanations.
On the one hand, many of our instances admit an optimal solution that uses two or three bins.
That is, a few Benders cuts are sufficient to achieve a matching dual bound.
On the other hand, we observed that the time needed by \setting{EFcons} to find an optimal solution is much higher than for \setting{pool}.
We explain the latter by the very expensive separation problem of Benders cuts, which requires to solve a sphere packing problem.
Since, as argued above, \setting{EFcons} will never encounter symmetric infeasible solutions, every separation problem is very time consuming.
Setting \setting{EFrow} and \setting{pool}, however, can explore symmetric parts of the branch-and-bound tree, which arguably increases the chance that heuristics find good feasible solution, while symmetric solutions can be separated very easily by making use of a symmetric cut pool.
Table~\ref{tab:statistics} confirms that the mean time needed to separate a cut using \setting{pool} and \setting{EFrow} is much smaller than for \setting{EFcons} when comparing the time spent for cut separation and the number of separated solutions.

In the context of bin packing, we can observe in Table~\ref{tab:conflict-cuts-statistics} that lifting the conflict cuts improves the solving time in almost all settings.
Again, the difference is more pronounced with rectangle packing than in sphere packing.
The largest gap appears to be between \setting{pool} and \setting{lifted pool} with up to almost~\SI{50}{\percent} time reduction for~\rectanglepack, and between the default conflict cuts and \setting{lifted} for~\spherepack, with up to \SI{20}{\percent} improvement. 
One unexpected trend we can observe in this table is that the improvement gained from lifting the cuts does not surpass the benefit of the compact extended formulation.
In particular, the gap between \setting{NC EFcons} and \setting{lifted EFcons} is smaller than with \setting{EFcons}.
One possible reason for this difference could come from the sheer difference in variables in the compact cuts.
In our test sets, there can be between $10$ to $25$ variables per non-compact cut~\eqref{eq:non-compact-vs-compact} after lifting, whereas the compact versions have at most five variables.

To answer~\eqref{Q2}, the results show that problems in which a small number of Benders cuts suffices to derive a good dual bound benefit from \setting{pool}; problems that require a lot of Benders cuts to prove optimality benefit from \setting{EFcons} (e.g., \setting{EFcons} reduces the number of separated solutions for \scheduling from more than two millions to less than~200).
This suggests that a hybrid strategy could be beneficial that first uses the \setting{pool} approach to explore the branch-and-bound tree to find good solutions, and then switches to \setting{EFcons} to keep the tree small.

\begin{table}[t]
  \caption{Statistics on different approaches for solving problems with conflict cuts.}
  \label{tab:conflict-cuts-statistics}
  \centering
  \begin{scriptsize}
  \begin{tabular*}{\textwidth}{@{}l@{\;\;\extracolsep{\fill}}rrrrrrrr@{}}
    \toprule
    & \multicolumn{4}{c}{with symmetry handling} & \multicolumn{4}{c}{without symmetry handling}\\
    \cmidrule{2-5} \cmidrule{6-9}
    setting & \#solved & time & cut-time & \#sepa & \#solved & time & cut-time & \#sepa \\
    \midrule
    \multicolumn{9}{@{}l}{\rectanglepack:}\\
    plain & \num{ 24} & \num{2678.0} & \num{2613.6} & \num{     894.4} & \num{ 24} & \num{2704.8} & \num{2642.1} & \num{    1022.0} \\
    lifted & \num{ 31} & \num{2340.6} & \num{2286.5} & \num{     858.4} & \num{ 24} & \num{2776.6} & \num{2693.3} & \num{    1046.3} \\
    pool & \num{ 35} & \num{1159.9} & \num{ 395.0} & \num{    2600.1} & \num{ 26} & \num{1960.4} & \num{ 498.5} & \num{    6981.4} \\
    lifted pool & \num{ 46} & \num{ 577.3} & \num{ 266.8} & \num{     911.4} & \num{ 32} & \num{1278.6} & \num{ 389.4} & \num{    2621.1} \\
    NC EFcons & \num{ 54} & \num{ 303.5} & \num{ 138.5} & \num{     119.6} & \num{ 36} & \num{1183.5} & \num{ 310.8} & \num{     254.3} \\
    lifted EFcons & \num{ 53} & \num{ 295.5} & \num{ 124.9} & \num{     123.6} & \num{ 37} & \num{1173.8} & \num{ 281.3} & \num{     235.2} \\
    EFcons & \num{ 59} & \num{ 232.9} & \num{ 118.4} & \num{     103.7} & \num{ 41} & \num{ 793.3} & \num{ 233.7} & \num{     199.0} \\
    \midrule
    \multicolumn{9}{@{}l}{\spherepack:}\\
    plain & \num{ 46} & \num{ 287.4} & \num{ 287.1} & \num{       8.8} & \num{ 42} & \num{ 357.9} & \num{ 357.8} & \num{      11.6} \\
    lifted & \num{ 50} & \num{ 242.8} & \num{ 242.5} & \num{       7.5} & \num{ 44} & \num{ 287.8} & \num{ 287.7} & \num{       9.3} \\
    pool & \num{ 57} & \num{ 109.2} & \num{ 109.0} & \num{      10.1} & \num{ 55} & \num{ 127.8} & \num{ 113.3} & \num{      24.3} \\
    lifted pool & \num{ 56} & \num{ 102.3} & \num{ 102.2} & \num{       7.2} & \num{ 56} & \num{ 104.6} & \num{ 103.0} & \num{      10.3} \\
    NC EFcons & \num{ 55} & \num{ 171.6} & \num{ 171.2} & \num{       5.7} & \num{ 54} & \num{ 198.3} & \num{ 197.7} & \num{       6.7} \\
    lifted EFcons & \num{ 55} & \num{ 168.8} & \num{ 168.5} & \num{       5.5} & \num{ 54} & \num{ 195.5} & \num{ 194.7} & \num{       6.9} \\
    EFcons & \num{ 55} & \num{ 173.4} & \num{ 173.1} & \num{       5.7} & \num{ 53} & \num{ 187.4} & \num{ 187.0} & \num{       6.6} \\
    \bottomrule
  \end{tabular*}
  \end{scriptsize}
\end{table}

\subsection*{Conclusion}
We have addressed the topic of detecting and exploiting symmetries in BD.
While automatic symmetry handling became a powerful component of modern MIP solvers~\cite{PfetschRehn2019}, BD cannot benefit from these techniques because MIP solvers are not aware of the symmetries of the cuts that are separated.
We demonstrated that symmetry detection for BD is, in principle, possible by creating an SDG that captures the symmetries of Benders cuts, and that it can be easily implemented within the solver \scip.
This allowed us to leverage state-of-the-art symmetry handling methods for MIP to BD to achieve substantial performance improvements.
Moreover, we systematically compared (to the best of our knowledge) for the first time the cut pool and extended formulation approach to enhance the separation of symmetric Benders cuts.
Based on our numerical study, we have seen that neither approach is dominant and suggested a hybrid approach that aims to create synergies between these approaches.

These findings provide fundamental insights into how symmetries should be exploited in BD.
This is particularly relevant for automatic BD schemes for generating Benders cuts that exist, e.g., for the solvers \solver{Cplex}~\cite{BonamiEtAl2020}, \solver{SAS}~\cite{SAS}, and \scip~\cite{Maher2021}.
Our long term goal is to incorporate our findings into such an automatic BD scheme, including a fully automatic symmetry detection algorithm for generic BD.
%
%
\bibliographystyle{plainurl}
\bibliography{bib}

@article{adulyasakEtAl2015,
  author =	 {Adulyasak, Yossiri and Cordeau, Jean-Fran\c{c}ois and Jans, Raf},
  title =	 {Benders Decomposition for Production Routing Under Demand Uncertainty},
  journal =	 {Operations Research},
  volume =	 {63},
  number =	 {4},
  pages =	 {851--867},
  year =	 {2015}
}

@article{akpinar2017combinatorial,
  title={Combinatorial {Benders} cuts for assembly line balancing problems with setups},
  author={Akpinar, Sener and Elmi, Atabak and Bekta{\c{s}}, Tolga},
  journal={European Journal of Operational Research},
  volume={259},
  number={2},
  pages={527--537},
  year={2017},
  publisher={Elsevier}
}

@inproceedings{AndersEtAl2023,
  author = {Anders, Markus and Schweitzer, Pascal and Stie{\ss}, Julian},
  title = {Engineering a Preprocessor for Symmetry Detection},
  booktitle = {21st International Symposium on Experimental Algorithms (SEA 2023)},
  series = {LIPIcs},
  year = {2023},
  doi = {10.4230/LIPIcs.SEA.2023.1}
}

@article{balas1978facets,
  title={Facets of the knapsack polytope from minimal covers},
  author={Balas, Egon and Zemel, Eitan},
  journal={SIAM Journal on Applied Mathematics},
  volume={34},
  number={1},
  pages={119--148},
  year={1978},
  publisher={SIAM}
}

@Article{benders1962,
  author =	 {Jacques F. Benders},
  title =	 {Partitioning procedures for solving mixed-variables programming problems},
  journal =	 {Numerische Mathematik},
  year =	 {1962},
  volume =	 {4},
  pages =	 {238--252}
 }

@article{BendottiEtAl2021,
  Title =	 {Orbitopal fixing for the full \mbox{(sub-)} orbitope and application to the Unit
                  Commitment Problem},
  Author =	 {Pascale Bendotti and Pierre Fouilhoux and C\'ecile Rottner},
  journal =	 {Mathematical Programming},
  volume =	 {186},
  year =	 {2021},
  pages =	 {337--372},
  doi =		 {10.1007/s10107-019-01457-1}
}

@unpublished{bendotti2026two,
  author = {Bendotti, Pascale and Heintzmann, Alexandre and Rottner, Cécile},
  title  = {A Two-Phase Branch\&Cut Scheme for Integer Programs with Polyhedral Symmetries},
  year   = {2026},
  note   = {Preprint submitted to HAL, hal-05514821v2},
  url    = {https://hal.science/hal-05514821}
}

@InProceedings{BonamiEtAl2020,
  author =	 "Bonami, Pierre and Salvagnin, Domenico and Tramontani, Andrea",
  editor =	 "Bienstock, Daniel and Zambelli, Giacomo",
  title =	 "Implementing Automatic {Benders} Decomposition in a Modern {MIP} Solver",
  booktitle =	 "Integer Programming and Combinatorial Optimization",
  year =	 "2020",
  pages =	 "78--90",
  isbn =	 "978-3-030-45771-6"
}

@article{buchheim2005linear,
  title={Linear optimization over permutation groups},
  author={Buchheim, Christoph and J{\"u}nger, Michael},
  journal={Discrete Optimization},
  volume={2},
  number={4},
  pages={308--319},
  year={2005},
  publisher={Elsevier}
}

@article{Daryalal2023,
  author =	 {Daryalal, Maryam and Pouya, Hamed and DeSantis, Marc Antoine},
  title =	 {Network Migration Problem: A Hybrid Logic-Based {Benders} Decomposition Approach},
  journal =	 {INFORMS Journal on Computing},
  volume =	 {35},
  number =	 {3},
  pages =	 {593--613},
  year =	 {2023},
  doi =		 {10.1287/ijoc.2023.1280},
}

@article{delorme2017,
  title =	 {Logic based {Benders'} decomposition for orthogonal stock cutting problems},
  journal =	 {Computers \& Operations Research},
  volume =	 {78},
  pages =	 {290--298},
  year =	 {2017},
  issn =	 {0305-0548},
  doi =		 {10.1016/j.cor.2016.09.009},
  author =	 {Maxence Delorme and Manuel Iori and Silvano Martello},
}

@article{DoornmalenHojny2024a,
  author =	 {van Doornmalen, Jasper and Hojny, Christopher},
  title =	 {A unified framework for symmetry handling},
  journal =	 {Mathematical Programming},
  volume =	 {212},
  pages =	 {217--271},
  year =	 {2025},
  doi = {10.1007/s10107-024-02102-2}
}

@article{fazel2012using,
  title={Using logic-based {Benders} decomposition to solve the capacity- and distance-constrained plant location problem},
  author={Fazel-Zarandi, Mohammad M. and Beck, J. Christopher},
  journal={INFORMS Journal on Computing},
  volume={24},
  number={3},
  pages={387--398},
  year={2012},
  publisher={INFORMS}
}

@article{guoEtAl2021,
  author =	 {Guo, Cheng and Bodur, Merve and Aleman, Dionne M. and Urbach, David R.},
  title =	 {Logic-Based {Benders} Decomposition and Binary Decision Diagram Based Approaches for
                  Stochastic Distributed Operating Room Scheduling},
  journal =	 {INFORMS Journal on Computing},
  volume =	 {33},
  number =	 {4},
  pages =	 {1551--1569},
  year =	 {2021},
  doi =		 {10.1287/ijoc.2020.1036},
}

@article{Han2025,
  title =	 {Integrated optimization of train makeup problem and resource scheduling in railway
                  marshalling yards: A hybrid {MILP-CP} approach with Logic-based {Benders}
                  decomposition},
  journal =	 {Transportation Research Part B: Methodological},
  volume =	 {200},
  pages =	 {103306},
  year =	 {2025},
  doi =    {10.1016/j.trb.2025.103306},
  author =	 {Peiran Han and Lingyun Meng and Xiaojie Luan and Nikola Bešinović and Jianrui Miao
                  and Yihui Wang and Zhengwen Liao},
}

@article{heching2019logic,
  title={A logic-based {Benders} approach to home healthcare delivery},
  author={Heching, Aliza and Hooker, John N. and Kimura, Ryo},
  journal={Transportation Science},
  volume={53},
  number={2},
  pages={510--522},
  year={2019},
  publisher={INFORMS}
}

@article{hojny2020knapsack,
  title={Knapsack polytopes: a survey},
  author={Hojny, Christopher and Gally, Tristan and Habeck, Oliver and L{\"u}then, Hendrik and Matter, Frederic and Pfetsch, Marc E. and Schmitt, Andreas},
  journal={Annals of Operations Research},
  volume={292},
  pages={469--517},
  year={2020},
  publisher={Springer}
}

@article{hojny2025detecting,
  title =	 {Detecting and handling reflection symmetries in mixed-integer (nonlinear)
                  programming and beyond},
  author =	 {Hojny, Christopher},
  journal =	 {Mathematical Programming Computation},
  pages =	 {31--78},
  year =	 {2026},
  volume =	 {18},
  doi =		 {10.1007/s12532-025-00289-9}
}

@InProceedings{hojny2026IPCO,
  author="Hojny, Christopher and Roy, C{\'e}dric",
  editor="Dey, Santanu S. and Di Summa, Marco and Salvagnin, Domenico",
  title="A Framework for Handling and Exploiting Symmetry in {Benders} Decomposition",
  booktitle="Integer Programming and Combinatorial Optimization",
  year="2026",
  publisher="Springer Nature Switzerland",
  address="Cham",
  pages="441--456",
  isbn="978-3-032-28691-8"
}

@article{hooker2003logic,
  title={Logic-based {Benders} decomposition},
  author={Hooker, John N. and Ottosson, Greger},
  journal={Mathematical Programming},
  volume={96},
  number={1},
  pages={33--60},
  year={2003},
  publisher={Springer}
}

@Book{hooker2024,
  author = 	 {John Hooker},
  title = 	 {Logic-Based {Benders} Decomposition},
  publisher = 	 {Springer},
  year = 	 {2024},
  series = 	 {Synthesis Lectures on Operations Research and Applications},
  address = 	 {Cham},
  edition = 	 {1}
}

@article{KaibelEtAl2011,
  title                    = {Orbitopal fixing},
  journal                  = {Discrete Optimization},
  volume                   = {8},
  number                   = {4},
  pages                    = {595--610},
  year                     = {2011},
  issn                     = {1572-5286},
  doi                      = {10.1016/j.disopt.2011.07.001},
  author                   = {Volker Kaibel and Matthias Peinhardt and Marc E. Pfetsch}
}

@article{kaparis2010separation,
  title={Separation algorithms for 0-1 knapsack polytopes},
  author={Kaparis, Konstantinos and Letchford, Adam N.},
  journal={Mathematical Programming},
  volume={124},
  pages={69--91},
  year={2010},
  publisher={Springer}
}

@article{khajavirad2024circle,
  title={The circle packing problem: A theoretical comparison of various convexification techniques},
  author={Khajavirad, Aida},
  journal={Operations Research Letters},
  volume={57},
  pages={107197},
  year={2024},
  publisher={Elsevier}
}

@article{lehmann2023accelerated,
  title={Accelerated {Benders} Decomposition for Variable-Height Transport Packaging Optimisation},
  author={Lehmann, Alain and Kleiminger, Wilhelm and Invernizzi, Hakim and Gautschi, Aurel},
  journal={arXiv preprint arXiv:2308.01104},
  year={2023}
}

@article{liberti2012reformulations,
  title={Reformulations in mathematical programming: automatic symmetry detection and exploitation},
  author={Liberti, Leo},
  journal={Mathematical Programming},
  volume={131},
  number={1},
  pages={273--304},
  year={2012},
  publisher={Springer}
}

@article{LibertiOstrowski2014,
  author  = {Leo Liberti and James Ostrowski},
  title   = {Stabilizer-based symmetry breaking constraints for mathematical programs},
  journal = {Journal of Global Optimization},
  year    = {2014},
  volume  = {60},
  pages   = {183--194}
}

@article{Maher2021,
  title =	 {Implementing the branch-and-cut approach for a general purpose {Benders’}
                  decomposition framework},
  journal =	 {European Journal of Operational Research},
  volume =	 {290},
  number =	 {2},
  pages =	 {479--498},
  year =	 {2021},
  issn =	 {0377-2217},
  doi =		 {10.1016/j.ejor.2020.08.037},
  author =	 {Stephen J. Maher}
}

@article{margot2009symmetry,
  title={Symmetry in integer linear programming},
  author={Margot, Fran{\c{c}}ois},
  journal={50 years of integer programming 1958--2008: from the early years to the state-of-the-art},
  pages={647--686},
  year={2009},
  publisher={Springer}
}

@article{Nauty,
  title =	 {Practical graph isomorphism, {II}},
  journal =	 {Journal of Symbolic Computation},
  volume =	 {60},
  pages =	 {94--112},
  year =	 {2014},
  doi =		 {10.1016/j.jsc.2013.09.003},
  author =	 {Brendan D. McKay and Adolfo Piperno}
}

@article{nemhauser1994lifted,
  title={Lifted cover facets of the 0--1 knapsack polytope with {GUB} constraints},
  author={Nemhauser, George L and Vance, Pamela H},
  journal={Operations Research Letters},
  volume={16},
  number={5},
  pages={255--263},
  year={1994},
  publisher={Elsevier}
}

@article{OstrowskiEtAl2011,
  year                     = {2011},
  issn                     = {0025-5610},
  journal                  = {Mathematical Programming},
  volume                   = {126},
  number                   = {1},
  doi                      = {10.1007/s10107-009-0273-x},
  title                    = {Orbital branching},
  author                   = {Ostrowski, James and Linderoth, Jeff and Rossi, Fabrizio and Smriglio, Stefano},
  pages                    = {147--178}
}

@article{padberg1975note,
  title={A note on zero-one programming},
  author={Padberg, Manfred W.},
  journal={Operations Research},
  volume={23},
  number={4},
  pages={833--837},
  year={1975},
  publisher={INFORMS}
}

@article{PfetschRehn2019,
  doi =		 {10.1007/s12532-018-0140-y},
  year =	 2019,
  volume =	 {11},
  number =	 {1},
  pages =	 {37--93},
  author =	 {Pfetsch, Marc E. and Rehn, Thomas},
  title =	 {A computational comparison of symmetry handling methods for mixed integer
                  programs},
  journal =	 {Mathematical Programming Computation}
}

@article{pisinger2007using,
  title={Using decomposition techniques and constraint programming for solving the two-dimensional bin-packing problem},
  author={Pisinger, David and Sigurd, Mikkel},
  journal={INFORMS Journal on Computing},
  volume={19},
  number={1},
  pages={36--51},
  year={2007},
  publisher={INFORMS}
}

@misc{releaseBenders,
  author = {Hojny, Christopher and Roy, C\'edric},
  title = {Supplementary material for the article ``{A} Framework for Handling and Exploiting Symmetry in {Benders} Decomposition''}, 
  year = {2026},
  doi = {10.5281/zenodo.23057507}
}

@article{riise2016recursive,
  title={Recursive logic-based {Benders’} decomposition for multi-mode outpatient scheduling},
  author={Riise, Atle and Mannino, Carlo and Lamorgese, Leonardo},
  journal={European Journal of Operational Research},
  volume={255},
  number={3},
  pages={719--728},
  year={2016},
  publisher={Elsevier}
}

@article{Salvagnin2005,
  title   = {A dominance procedure for integer programming},
  author  = {Salvagnin, Domenico},
  journal = {Master’s thesis, University of Padova, Padova, Italy},
  year    = {2005}
}

@inproceedings{Salvagnin2018,
  doi =		 {10.1007/978-3-319-93031-2_37},
  year =	 {2018},
  publisher =	 {Springer International Publishing},
  pages =	 {521--529},
  author =	 {Domenico Salvagnin},
  title =	 {Symmetry Breaking Inequalities from the {Schreier-Sims} Table},
  booktitle =	 {Integration of Constraint Programming, Artificial Intelligence, and Operations Research},
  editor =	 {van Hoeve, Willem-Jan},
  isbn =	 {978-3-319-93031-2}
}

@Manual{SAS,
  title = 	 {{SAS} {Programming} {Documentation}},
  organization = {SAS Institute Inc.},
  note = 	 {\url{https://documentation.sas.com/doc/en/pgmsascdc/v_068/casmopt/casmopt_benders_overview.htm}},
}

@article{SCIP8,
  author = {Bestuzheva, Ksenia and Besan\c{c}on, Mathieu and Chen, Wei-Kun and Chmiela, Antonia and Donkiewicz, Tim and van Doornmalen, Jasper and Eifler, Leon and Gaul, Oliver and Gamrath, Gerald and Gleixner, Ambros and Gottwald, Leona and Graczyk, Christoph and Halbig, Katrin and Hoen, Alexander and Hojny, Christopher and van der Hulst, Rolf and Koch, Thorsten and L\"{u}bbecke, Marco and Maher, Stephen J. and Matter, Frederic and M\"{u}hmer, Erik and M\"{u}ller, Benjamin and Pfetsch, Marc E. and Rehfeldt, Daniel and Schlein, Steffan and Schl\"{o}sser, Franziska and Serrano, Felipe and Shinano, Yuji and Sofranac, Boro and Turner, Mark and Vigerske, Stefan and Wegscheider, Fabian and Wellner, Philipp and Weninger, Dieter and Witzig, Jakob},
  title = {Enabling Research through the {SCIP Optimization Suite 8.0}},
  year = {2023},
  issue_date = {June 2023},
  publisher = {Association for Computing Machinery},
  address = {New York, NY, USA},
  volume = {49},
  number = {2},
  issn = {0098-3500},
  url = {https://doi.org/10.1145/3585516},
  doi = {10.1145/3585516},
  journal = {ACM Trans. Math. Softw.},
  articleno = {22},
  numpages = {21}
}

@misc{SCIP9,
  title =	 {The {SCIP} {Optimization} {Suite} 9.0},
  author =	 {Suresh Bolusani and Mathieu Besançon and Ksenia Bestuzheva and Antonia Chmiela and
                  João Dionísio and Tim Donkiewicz and Jasper van Doornmalen and Leon Eifler and
                  Mohammed Ghannam and Ambros Gleixner and Christoph Graczyk and Katrin Halbig and
                  Ivo Hedtke and Alexander Hoen and Christopher Hojny and Rolf van der Hulst and
                  Dominik Kamp and Thorsten Koch and Kevin Kofler and Jurgen Lentz and Julian Manns
                  and Gioni Mexi and Erik Mühmer and Marc E. Pfetsch and Franziska Schlösser and
                  Felipe Serrano and Yuji Shinano and Mark Turner and Stefan Vigerske and Dieter
                  Weninger and Lixing Xu},
  year =	 {2024},
  eprint =	 {2402.17702},
  archivePrefix ={arXiv},
  primaryClass = {math.OC}
}

@article{tang1988models,
  title={Models arising from a flexible manufacturing machine, {Part} {I}: minimization of the number of tool switches},
  author={Tang, Christopher S. and Denardo, Eric V.},
  journal={Operations Research},
  volume={36},
  number={5},
  pages={767--777},
  year={1988},
  publisher={INFORMS}
}

@article{tran2016decomposition,
  title={Decomposition methods for the parallel machine scheduling problem with setups},
  author={Tran, Tony T. and Araujo, Arthur and Beck, J. Christopher},
  journal={INFORMS Journal on Computing},
  volume={28},
  number={1},
  pages={83--95},
  year={2016},
  publisher={INFORMS}
}

@article{van2023traditional,
  title={A traditional {Benders’} approach to sports timetabling},
  author={Van Bulck, David and Goossens, Dries},
  journal={European Journal of Operational Research},
  volume={307},
  number={2},
  pages={813--826},
  year={2023},
  publisher={Elsevier}
}

@Article{WachterBiegler2006,
  author = 	 {Andreas W\"achter and Lorenz T. Biegler},
  title = 	 {On the implementation of an interior-point filter line-search algorithm for large-scale nonlinear programming},
  journal = 	 {Mathematical Programming},
  year = 	 {2006},
  volume = 	 {106},
  pages = 	 {25--57},
  doi =		 {10.1007/s10107-004-0559-y}
}

@article{zemel1989easily,
  title={Easily computable facets of the knapsack polytope},
  author={Zemel, Eitan},
  journal={Mathematics of Operations Research},
  volume={14},
  number={4},
  pages={760--764},
  year={1989},
  publisher={INFORMS}
}

%
%

\newpage
\appendix

\section{Full Models for Tested Benders Decompositions}\label{appendix:MIPS}

This appendix provides the full details for the master problems and subproblems used in our implementation of BD for the rectangle packing and scheduling problem.
The presented models are also the models that we have used in our implementation.

\subsection{Rectangle Packing}\label{app:rectangle-packing}
In the rectangle packing problem, we are given a set~$N$ of rectangles and a set~$B$ of rectangular bins.
Rectangle~$j \in N$ has width~$w_j$ and height~$h_j$, whereas bin~$i \in B$ has width~$W_i$ and height~$H_i $.
The task is to find an assignment of the rectangles~$N$ to the bins~$B$ such that the rectangles assigned to the same bin~$i$ can be placed inside bin~$i$ without overlap.
Among all such assignments, one is looking for one that minimizes the number of bins with an assigned rectangle.
We furthermore operate under the model where all rectangles have edges horizontally or vertically aligned and cannot be rotated.

The master problem follows the blueprint problem~\eqref{eq:naive-MK} and introduces the usual binary variables~$x \in \B{B \times N}$ to indicate whether rectangle~$j \in N$ is assigned to bin~$i$.
Additional variables~$u_i \in \B{}$, $i \in B$, are introduced to indicate whether bin~$i$ is non-empty.
Instead of deciding feasibility of an assignment purely based on solving Benders subproblems, we enhance the master problem by ruling out assignments of rectangles whose total area exceeds the area of the bin.
The master problem is then given by
\begin{subequations}
  \label{eq:master-2DRP}
  \begin{align}
    \min \sum_{i \in B} u_i &&& \nonumber\\
    \sum_{i \in B} x_{ij} &= 1, && j \in N,\\
    \sum_{j \in C} x_{ij} &\leq \card{C} - 1, && (C,i) \in \mathcal{I},\\
    \sum_{j \in N} (w_j \cdot h_j) x_{ij} &\leq W_i\cdot H_i, && i \in B,\\
    x_{ij} &\leq u_i, && i \in B,\; j \in N,\\
    x_{ij} &\in \B{}, && i \in B,\; j \in N,\\
    u_i &\in \B{}, && i \in B.
  \end{align}
\end{subequations}

The family of Benders cuts~$\sum_{j \in C} x_{ij} \leq \card{C} - 1$ for~$(C,i) \in \mathcal{I}$ can then be separated by taking a solution~$\fixed{x} \in \B{B \times N}$ and checking whether the items assigned to bin~$i \in B$ fit into the bin.
Following~\cite{pisinger2007using}, given a bin~$i$ and solution~$\fixed{x}$, the subproblem for deciding feasibility of~$\fixed{x}$ can be formulated as
\begin{subequations}
  \label{eq:subproblem-2DRP}
  \begin{align}
    \max 0 &&& \nonumber \\
    y^{(1)}_j + \fixed{x}_{ij}w_j + l_{jk} W_i &\leq y^{(1)}_k + W_i, && j\neq k\in N, \label{eq:sub-box-hspace} \\ 
    y^{(2)}_j + \fixed{x}_{ij}h_j + b_{jk} H_i &\leq y^{(2)}_k + H_i, && j\neq k\in N, \label{eq:sub-box-vspace} \\
    l_{jk} + l_{kj} + b_{jk} + b_{kj} &\geq 1, && j < k\in N, \label{eq:sub-box-non-overlap} \\
    l_{jk} + l_{kj} &\leq 1, && j < k\in N \label{eq:sub-box-left-right}\\
    b_{jk} + b_{kj} &\leq 1, && j < k\in N \label{eq:sub-box-top-bottom} \\
    l_{jk}, b_{jk} &\in \{0, 1\}, && j\neq k \in N, \\
    0 \leq y^{(1)}_j &\leq \fixed{x}_{ij}(W_i - w_j), && j\in N, \\
    0 \leq y^{(2)}_j &\leq \fixed{x}_{ij}(H_i - h_j), && j\in N,
  \end{align}
\end{subequations}
where the pair~$(y^{(1)}_j, y^{(2)}_j)$ acts as the bottom left corner of the~$j$-th rectangle.
Observe that if~$\fixed{x}_{ij} = 0$ for some~$j \in N$, then the rectangle~$j$ is forced to be placed at~$(0, 0)$ and its width and height are ignored. 
In practice, one would remove the rectangle~$j$ from the subproblem.
For the sake of exposition, however, we decided to keep all~$x$-variables in the subproblem because it allows to define an SDG for the subproblem by identifying~$\hat{x}$ with the corresponding master variables.

The binary variable~$l_{jk}$ indicates whether rectangles~$j$ and~$k$ do not overlap horizontally, and~$j$ is to the left of~$k$.
Indeed, if~$l_{jk}=1$, Equation~\eqref{eq:sub-box-hspace} reduces to~$y^{(1)}_j + w_j \leq y^{(1)}_k$.
Otherwise, if~$l_{jk} =0$, the constraint is trivially always satisfied.
Analogously, the variable~$b_{jk}$ is 1 only if rectangles~$j$ and~$k$ do not overlap vertically, and~$j$ is below~$k$.
This means that for two rectangles to overlap, they necessarily would have~$l_{jk} = l_{kj} = b_{jk} = b_{kj} = 0$, which is prevented by Equation~\eqref{eq:sub-box-non-overlap}.
Moreover, we add~\eqref{eq:sub-box-left-right} and~\eqref{eq:sub-box-top-bottom} to the model to directly encode that neither rectangle~$j$ is simultaneously to the left and right of rectangle~$k$ nor above/below it.

To build an SDG for the \rectanglepack problem, we follow Theorem~\ref{thm:merge-SDGs} and create an SDG for~\eqref{eq:master-2DRP} and~\eqref{eq:subproblem-2DRP}.
Afterwards, these two SDGs are merged into a single SDG by identifying variable nodes with each other.
Since both problems are MIPs, these SDGs can be created following the standard rules for building SDGs for MIPs~\cite{Salvagnin2005}.
We can, however, create a more compact SDG~$G=(V, E)$ for the subproblem in a similar way to that of \spherepack:
We define the node set
\[
  V = \{a\} \cup \bigcup_{j \in N} \{v_j, u^w_j, u^h_j\},
\]
where~$a$ is the anchor of the graph, $v_j$ is a variable node representing variable~$x_{ij}$, and~$u^w_j$ and~$u^h_j$ are nodes representing the width and height of item~$j$, respectively.
The anchor~$a$ receives a color that uniquely characterizes the dimensions~$(W_i,H_i)$ of the respective bin~$i$, the variable nodes are colored according to their objective coefficient and type in the master problem, and the nodes~$u^w_j$ and~$u^h_j$ receive a unique color that indicates them as width-node and height-node, respectively.
Moreover, we define the edge set
\[
  E = \bigcup_{j \in N} \brackets{
    \brackets{a,u^w_j}, \brackets{u^w_j,u^h_j},\brackets{u^h_j,v_j}
    },
\]
where~$\{a,u^w_j\}$ receives a color corresponding to the width~$w_j$, $\{u^w_j,u^h_j\}$ receives a color based on the height~$h_j$, and edge~$\{u^h_j,v_j\}$ remains uncolored.
Due to this construction, an automorphism of~$G$ can only exchange variable nodes~$v_j$ and~$v_{j'}$ if~$w_j = w_{j'}$ and~$h_j = h_{j'}$.
Thus, the SDG~$G$ captures the symmetries of the subproblem.

Determining the value of the lifting coefficients in \rectanglepack is too costly as it is equivalent to solving a smaller instance of \rectanglepack.
We therefore restrict ourself to a smaller value of~$\lift'$, using similar observations as for \spherepack. 
It is easy to verify that~$\lift'_{ij} \geq 1$ is equivalent to~$(h_j \geq h_k)\wedge (w_j \geq w_k)$ for all~$k\in C$.
Checking if any two rectangles~$k_1, k_2$ fit inside rectangle~$j$ is equivalent to the disjunction~$(h_{k_1} + h_{k_2} \leq h_j) \lor (w_{k_1} + w_{k_2} \leq w_j)$ (assuming~$\lift'_{ij} \geq 1$ to begin with).
We can go even one step further.

\begin{lemma}
  \label{lemma:3fit}
  Let~$a, b, c$ be three non-rotatable rectangles.
  If they can be packed inside of a rectangle of dimensions~$H\times W$, the packing admits a vertical or horizontal guillotine cut.
  In other words, there is a reordering~$\brackets{i, j, k} = \brackets{a, b, c}$ such that boxes~$i$ and~$j$ fit inside a rectangle of dimension~$(H-h_k)\times W$ or~$H \times (W-w_k)$.
\end{lemma}
\begin{proof}
  For any given valid packing of~$\brackets{a, b, c}$, let~$G_H = (V, E_H)$ be the intersection graph for the projection of the packing on the height of the box. 
  In particular, $V=\brackets{a, b, c}$ and~$\brackets{i, j}\in E_H$ if the rectangles~$i$ and~$j$ overlap on that projection.
  Similarly, we construct~$G_W = (V, E_W)$ for the projection on the width of box.
  Observe that because the projections are onto a single dimension, the packing admits a vertical (resp.\ horizontal) guillotine cut if and only if~$G_W$ (resp.\ $G_H$) is not connected, see Figure~\ref{fig:example-guillotine} for an illustration.
  Additionally, because the packing has no overlap, $E_H \cap E_W = \emptyset$.
  In other words, $|E_H| + |E_W| \leq 3$.
  As~$G_H$ and~$G_W$ both necessitate at least two edges to be connected, this means the packing must admit a guillotine cut.
\end{proof}
\begin{figure}
  \centering
  \begin{tikzpicture}[scale=0.5]
    \draw[thick] (0,-1) -- (11, -1);
    \draw[thick] (11, -1) -- (11, 10);
    \filldraw[fill=col5] (1, 9) rectangle node[circle, fill=col5!50]{$a$} (4, 4);
    \filldraw[fill=col3] (2, 3) rectangle node[circle, fill=col3!50]{$b$} (7, 0);
    \filldraw[fill=col4] (8, 1) rectangle node[circle, fill=col4!50]{$c$} (10, 7);

    \draw[decorate, decoration={calligraphic brace, amplitude=5pt}, thick] (1, -0.8) -- (4, -0.8) node[pos=0.5, above=2pt]{$a$};
    \draw[decorate, decoration={calligraphic brace, amplitude=5pt}, thick] (7, -1.2) -- (2, -1.2) node[pos=0.5, below=2pt]{$b$};
    \draw[decorate, decoration={calligraphic brace, amplitude=5pt}, thick] (8, -0.8) -- (10, -0.8) node[pos=0.5, above=2pt]{$c$};
    \draw[decorate, decoration={calligraphic brace, amplitude=5pt}, thick] (11.2, 9) -- (11.2, 4) node[pos=0.5, right=2pt]{$a$};
    \draw[decorate, decoration={calligraphic brace, amplitude=5pt}, thick] (11.2, 3) -- (11.2, 0) node[pos=0.5, right=2pt]{$b$};
    \draw[decorate, decoration={calligraphic brace, amplitude=5pt}, thick] (10.8, 1) -- (10.8, 7) node[pos=0.5, left=2pt]{$c$};

    \begin{scope}[line width=3pt]
      \draw[col5]  (1, -1) -- ( 2, -1); 
      \draw[black] (2, -1) -- ( 4, -1);
      \draw[col3]  (4, -1) -- ( 7, -1);
      \draw[col4]  (8, -1) -- (10, -1);
      
      \draw[col5]  (11, 7) -- (11, 9);
      \draw[black] (11, 7) -- (11, 4);
      \draw[col4]  (11, 4) -- (11, 3);
      \draw[black] (11, 3) -- (11, 1);
      \draw[col3]  (11, 1) -- (11, 0);
    \end{scope}

    \node (a1) at (17, 7) [circle, draw=black] {};
    \node (b1) at (20, 7) [circle, draw=black] {};
    \node (c1) at (23, 7) [circle, draw=black] {};
    \node (a2) at (17, 3) [circle, draw=black] {};
    \node (b2) at (20, 3) [circle, draw=black] {};
    \node (c2) at (23, 3) [circle, draw=black] {};
    
    \draw (a1) node[ left=6pt]{$G_W$};
    \draw (a2) node[ left=6pt]{$G_H$};
    \draw (a1) node[above=4pt]{$a$};
    \draw (b1) node[above=4pt]{$b$};
    \draw (c1) node[above=4pt]{$c$};
    \draw (a2) node[below=4pt]{$a$};
    \draw (b2) node[below=4pt]{$b$};
    \draw (c2) node[below=4pt]{$c$};

    \draw (a1) -- (b1);
    \draw (a2) to[out=20,  in=160] (c2);
    \draw (b2) -- (c2);

    \begin{scope}[col2, line width=2pt, dashed]
      \draw ( 7.5, -1.7) -- ( 7.5, 9.5);
      \draw (21.5,  5.5) -- (21.5, 8.5);
    \end{scope}

  \end{tikzpicture}
  \label{fig:example-guillotine}
  \caption{Example of guillotine cut, in dashed orange (left) corresponding to a cut in the intersection graph of the horizontal projection (right).}
\end{figure}
We can therefore test if a given set of three rectangles fit inside a larger fourth one as it reduces to six tests of fitting two rectangles inside a larger third, which we can do in constant time.
In other words, we can efficiently compute~$\lift^3_{ij}$.

\subsection{Machine Scheduling with Setup Times}\label{app:MS}

Let~$M$ be a set of machines and~$N$ be a set of jobs.
For each job~$j \in N$, its processing time on machine~$i \in M$ is~$p_{ij}$.
Moreover, if job~$k$ is processed immediately after job~$j$ on machine~$i$, a setup time of~$s_{ijk}$ is needed.
The machine scheduling problem, requires to assign each job to exactly one machine.
The makespan of a machine is then determined by finding an ordering of the assigned jobs on that machine such that the total setup time is minimized; the makespan is then this minimum total setup time plus the total processing time of the assigned jobs.
The goal of the machine scheduling problem is to find an assignment of jobs to machines such that the maximum makespan of a machine is minimized.

\paragraph{Master Problem}
In our implementation of the master problem for the machine scheduling problem, we use an enhanced formulation as proposed by~\cite{tran2016decomposition}.
To this end, recall from Section~\ref{sec:Applications} that finding the optimal setup time of the jobs assigned to a machine requires to solve a TSP instance.
The idea of~\cite{tran2016decomposition} is therefore to embed a relaxation of the TSP into the master problem, which allows to encode a lower bound on the optimal setup time into the model.
To this end, the set of jobs~$N$ is augmented by an artificial job~$0$ that serves as the starting and finishing job on each machine;
denote~$N_0 = N \cup \{0\}$. We assume that both the processing time~$p_0$ and the setup times~$s_{ij0}$ from any job~$j \in N$ to job~$0$ is~0.
We do not impose~$s_{i0k} = 0$ though to allow for modeling start up times of machines.
Additionally, the setup times need to respect the triangle inequality, i.e., $s_{ijk} + s_{ikl} \geq s_{ijl}$ for all~$j, k, l\in N$.

The model of~\cite{tran2016decomposition} then introduces a variable~$T$ to model the maximum makespan.
Moreover, it has continuous variables~$y_{ijk}$ to model whether job~$j$ is the direct predecessor of job~$k$ on machine~$i$.
These variables can be considered as a relaxation of the variables of the subtour elimination formulation of the TSP on machine~$i$.
The model of~\cite{tran2016decomposition} then reads as follows:
\begin{subequations}
  \begin{align}
    \min T & \nonumber \\
    \sum_{i \in M} x_{ij} &= 1, && j \in N, \label{eq:master-MSPS-partition}\\
    x_{i0} &= 1, && i \in M, \label{eq:master-MSPS-job0}\\
    x_{ij} &= \sum_{k \in N_0} y_{ijk}, && j \in N_0,\; i \in M, \label{eq:master-MSPS-projleft}\\
    x_{ij} &= \sum_{j \in N_0} y_{ikj}, && j \in N_0,\; i \in M, \label{eq:master-MSPS-projright}\\
    \xi_i &= \sum_{j \in N_0, k \in N_0} y_{ijk} s_{ijk}, && i \in M, \label{eq:master-MSPS-setuptimes}\\
    \sum_{j \in N} x_{ij} p_{ij} + \xi_i &\leq T, && i \in M, \label{eq:master-MSPS-tottime}\\
    T_i(C) - \sum_{j\in C} (1-x_{ij})\theta_{ij} &\leq T, && (C, i)\in \mathcal{I}, \label{eq:master-benders}\\
    x_{ij}  &\in \{0,1\},     && j \in N,\; i \in M, \\
    y_{ijk} &\in [0,1],     && j, k \in N_0,\; i \in M,\\
    \xi_i   &\in \R_{\geq 0}, && i \in M.
  \end{align}
\end{subequations}
The variables~$x_{ij}$ model whether job~$j$ is assigned to machine~$i$, and~\eqref{eq:master-MSPS-partition} enforces that each job is assigned to exactly one machine.
Moreover, the auxiliary job~$0$ is assigned to every machine via~\eqref{eq:master-MSPS-job0}.
The~$y$-variables are linked with the~$x$-variables via~\eqref{eq:master-MSPS-projleft} and~\eqref{eq:master-MSPS-projright}.
These constraints enforce that, if a job~$j$ is assigned to machine~$i$, then the corresponding node in a relaxed TSP solution is entered and left exactly once (by the possibly fractional TSP solution~$y$).
The relaxed TSP solution~$y$ then allows to derive a lower bound on the optimal setup time on machine~$i$ via~\eqref{eq:master-MSPS-setuptimes}.
Constraint~\eqref{eq:master-MSPS-tottime} then links the lower bounds and the processing times of jobs assigned to machine~$i$, to derive a lower bound on the makespan~$T$.
Finally, \eqref{eq:master-benders} are the Benders cuts that provide the correct linking of the makespan variable~$T$ with the~$x$-variables.

\paragraph{Subproblem}
The Benders cuts are derived as follows.
Given a solution~$\fixed{x}$ of the master problem, the subproblem determines the minimum time required to process all jobs on the respective machines.
As the solution revolves around finding an optimal ordering, we can use a slight variation of the subtour elimination formulation of the TSP:
\begin{subequations}
  \label{eq:TSP}
  \begin{align}
    \min \sum_{j\in N_0} \fixed{x}_{ij} p_{ij} &+ \sum_{j\in N_0} \sum_{k\in N_0 \setminus\brackets{j}} s_{ijk} z_{jk} \nonumber\\
    \sum_{j\in N_0} z_{jk} &= 1, && k\in N_0, 
    \label{eq:TSP-in-degree}\\
    \sum_{k\in N_0} z_{jk} &= 1, && j\in N_0, \label{eq:TSP-out-degree}\\
    \sum_{j\in S} \sum_{k\in S\setminus\brackets{j}} z_{jk} &\leq |S| -1, && S\subsetneq \brackets{j\in N_0: \fixed{x}_{ij} = 1}, |S| \geq 2, \label{eq:TSP-subtour} \\
    z_{jj} &= 1 - \fixed{x}_{ij}, && j\in N_0, \label{eq:TSP-assignment} \\
    z_{jk} &\in \brackets{0, 1}, && j, k\in N_0.
  \end{align}
\end{subequations}
Formulation~\eqref{eq:TSP} represents each job as a node of a graph, and uses variable~$z_{jk}$ to indicate whether there is a directed edge from~$j$ to~$k$.
Equations~\eqref{eq:TSP-in-degree} and~\eqref{eq:TSP-out-degree} ensure that exactly one edge is coming in and going out of each node, and Constraint~\eqref{eq:TSP-subtour} guarantees that there cannot be subcycles on the set of selected jobs.
Because Constraint~\eqref{eq:TSP-assignment} forces the selected jobs to not be isolated, the variables~$z_{jk}$ necessarily form a cycle over the selected jobs.

\section{Detailed Numerical Results}
\label{app:detailedresults}

This appendix provides detailed numerical results for the settings \setting{pool}, \setting{EFrow}, and \setting{EFcons} (with symmetry handling enabled) for all instances of the test sets \rectanglepack (Table~\ref{tab:fulloptrp}), \spherepack (Table~\ref{tab:fulloptsp}), and \scheduling (Table~\ref{tab:fulloptsched}) that could be solved by all three settings within the time limit.
For each instance, these tables provide the running time in seconds (running time), the number of nodes in the branch-and-bound tree (\#nodes), and the number of times a solution was separated (\#sepa).
The last line of each table provides the shifted geometric mean of the respective columns.

\begin{table}[t]
  \caption{Detailed statistics for \rectanglepack instances that are solved by all Benders approaches.}
  \label{tab:fulloptrp}
  \centering
  \begin{scriptsize}
  \begin{tabular*}{\textwidth}{@{}l@{\;\;\extracolsep{\fill}}rrrrrrrrr@{}}
    \toprule
    & \multicolumn{3}{c}{running time} & \multicolumn{3}{c}{\#nodes} & \multicolumn{3}{c}{\#sepa}\\
    \cmidrule{2-4} \cmidrule{5-7} \cmidrule{8-10}
    instance & pool & EFrow & EFcons & pool & EFrow & EFcons & pool & EFrow & EFcons \\
    \midrule
  2DRP\_20N5grp\_0 & \num{1101.7} & \num{5577.6} & \num{ 304.4} & \num{      6104} & \num{     12814} & \num{      1575} & \num{      4993} & \num{      5081} & \num{       231}\\
  2DRP\_20N5grp\_1 & \num{   5.0} & \num{  80.2} & \num{ 747.6} & \num{         1} & \num{        39} & \num{        95} & \num{         5} & \num{        41} & \num{        76}\\
 2DRP\_20N5grp\_10 & \num{  36.5} & \num{  33.0} & \num{   6.1} & \num{     23620} & \num{       701} & \num{       200} & \num{     15254} & \num{       661} & \num{        93}\\
 2DRP\_20N5grp\_11 & \num{  13.1} & \num{  48.6} & \num{   9.9} & \num{      1552} & \num{      1617} & \num{       620} & \num{      1498} & \num{       985} & \num{       121}\\
 2DRP\_20N5grp\_13 & \num{   3.1} & \num{  48.6} & \num{  10.8} & \num{        25} & \num{       719} & \num{       128} & \num{        61} & \num{       589} & \num{       136}\\
 2DRP\_20N5grp\_17 & \num{  16.5} & \num{   7.8} & \num{   5.7} & \num{        93} & \num{        86} & \num{        52} & \num{       221} & \num{        67} & \num{        43}\\
 2DRP\_20N5grp\_18 & \num{  16.2} & \num{  39.4} & \num{  46.0} & \num{        28} & \num{        36} & \num{        46} & \num{        74} & \num{        51} & \num{        58}\\
 2DRP\_20N5grp\_19 & \num{ 113.9} & \num{   9.7} & \num{   6.6} & \num{        13} & \num{        47} & \num{        22} & \num{        57} & \num{        84} & \num{        48}\\
 2DRP\_20N5grp\_21 & \num{1662.7} & \num{6436.9} & \num{6446.2} & \num{        20} & \num{        11} & \num{        11} & \num{        57} & \num{        22} & \num{        22}\\
 2DRP\_20N5grp\_23 & \num{ 521.5} & \num{   2.5} & \num{   2.5} & \num{         9} & \num{         8} & \num{         8} & \num{        24} & \num{        12} & \num{        12}\\
 2DRP\_20N5grp\_28 & \num{  37.1} & \num{  52.0} & \num{  33.0} & \num{       508} & \num{       387} & \num{       264} & \num{       759} & \num{       211} & \num{       164}\\
 2DRP\_20N5grp\_29 & \num{  60.1} & \num{1088.2} & \num{  14.8} & \num{      1377} & \num{    633909} & \num{       375} & \num{      1886} & \num{      8383} & \num{       189}\\
  2DRP\_20N5grp\_3 & \num{3729.4} & \num{ 368.4} & \num{ 369.0} & \num{         1} & \num{         1} & \num{         1} & \num{         5} & \num{         7} & \num{         7}\\
 2DRP\_20N5grp\_30 & \num{ 148.5} & \num{ 101.6} & \num{ 101.3} & \num{         1} & \num{         1} & \num{         1} & \num{         5} & \num{         9} & \num{         9}\\
 2DRP\_20N5grp\_31 & \num{3180.0} & \num{   7.0} & \num{   7.0} & \num{         1} & \num{         1} & \num{         1} & \num{         5} & \num{         5} & \num{         5}\\
 2DRP\_20N5grp\_35 & \num{  40.8} & \num{  62.6} & \num{  62.9} & \num{         1} & \num{         1} & \num{         1} & \num{        18} & \num{         5} & \num{         5}\\
 2DRP\_20N5grp\_37 & \num{  17.8} & \num{ 242.5} & \num{  16.5} & \num{      1766} & \num{      2598} & \num{       206} & \num{      1802} & \num{      2230} & \num{       186}\\
  2DRP\_20N5grp\_5 & \num{  32.8} & \num{ 983.6} & \num{  44.2} & \num{       350} & \num{      7811} & \num{       264} & \num{       486} & \num{      4040} & \num{       135}\\
  2DRP\_20N5grp\_6 & \num{  29.3} & \num{3549.8} & \num{  28.3} & \num{       886} & \num{     40638} & \num{       388} & \num{       856} & \num{     14927} & \num{       157}\\
  2DRP\_20N5grp\_7 & \num{   9.1} & \num{  19.2} & \num{  19.2} & \num{         1} & \num{         1} & \num{         1} & \num{         5} & \num{         7} & \num{         7}\\
  2DRP\_20N5grp\_8 & \num{ 117.7} & \num{1204.0} & \num{  91.2} & \num{      6333} & \num{      6841} & \num{      5326} & \num{      7711} & \num{      3474} & \num{       225}\\
  2DRP\_20N5grp\_9 & \num{ 661.5} & \num{ 265.6} & \num{ 265.1} & \num{         3} & \num{         1} & \num{         1} & \num{        35} & \num{        10} & \num{        10}\\
  2DRP\_25N5grp\_0 & \num{ 423.3} & \num{1868.5} & \num{  47.8} & \num{      1270} & \num{      1373} & \num{       266} & \num{      1653} & \num{       800} & \num{       136}\\
  2DRP\_25N5grp\_1 & \num{  38.9} & \num{ 204.9} & \num{ 204.8} & \num{         1} & \num{        42} & \num{        42} & \num{        19} & \num{        15} & \num{        15}\\
 2DRP\_25N5grp\_17 & \num{  29.1} & \num{ 813.9} & \num{  45.7} & \num{       340} & \num{      6266} & \num{       547} & \num{       446} & \num{      3240} & \num{       145}\\
 2DRP\_25N5grp\_19 & \num{5975.0} & \num{  12.4} & \num{  13.2} & \num{        97} & \num{        34} & \num{        34} & \num{       196} & \num{        40} & \num{        40}\\
 2DRP\_25N5grp\_28 & \num{  13.2} & \num{ 223.8} & \num{  68.0} & \num{        57} & \num{       869} & \num{       781} & \num{       165} & \num{       627} & \num{       420}\\
 2DRP\_25N5grp\_31 & \num{  16.3} & \num{  84.9} & \num{  85.1} & \num{         1} & \num{         6} & \num{         6} & \num{        10} & \num{        11} & \num{        11}\\
 2DRP\_25N5grp\_35 & \num{  30.2} & \num{   1.8} & \num{   1.8} & \num{         1} & \num{         1} & \num{         1} & \num{        14} & \num{         7} & \num{         7}\\
\midrule
mean (29 inst.) & \num{  80.5} & \num{ 122.5} & \num{  41.3} & \num{      98.4} & \num{     221.9} & \num{      76.5} & \num{     131.5} & \num{     138.6} & \num{      46.5}\\
    \bottomrule
  \end{tabular*}
  \end{scriptsize}
\end{table}
\begin{table}[t]
  \caption{Detailed statistics for \spherepack instances that are solved by all Benders approaches.}
  \label{tab:fulloptsp}
  \centering
  \begin{scriptsize}
  \begin{tabular*}{\textwidth}{@{}l@{\;\;\extracolsep{\fill}}rrrrrrrrr@{}}
    \toprule
    & \multicolumn{3}{c}{running time} & \multicolumn{3}{c}{\#nodes} & \multicolumn{3}{c}{\#sepa}\\
    \cmidrule{2-4} \cmidrule{5-7} \cmidrule{8-10}
    instance & pool & EFrow & EFcons & pool & EFrow & EFcons & pool & EFrow & EFcons \\
    \midrule
  2DSP\_10N2grp\_0 & \num{ 112.9} & \num{ 357.3} & \num{ 358.1} & \num{         1} & \num{         1} & \num{         1} & \num{        16} & \num{         6} & \num{         6}\\
  2DSP\_10N2grp\_1 & \num{   1.3} & \num{   1.3} & \num{   1.3} & \num{         1} & \num{         1} & \num{         1} & \num{         1} & \num{         1} & \num{         1}\\
 2DSP\_10N2grp\_10 & \num{  12.2} & \num{  40.4} & \num{  40.2} & \num{         1} & \num{         5} & \num{         5} & \num{        15} & \num{        19} & \num{        19}\\
 2DSP\_10N2grp\_12 & \num{   0.9} & \num{   0.8} & \num{   0.9} & \num{         1} & \num{         1} & \num{         1} & \num{         1} & \num{         1} & \num{         1}\\
 2DSP\_10N2grp\_13 & \num{  16.4} & \num{  46.0} & \num{  49.4} & \num{         1} & \num{         1} & \num{         1} & \num{         8} & \num{         6} & \num{         6}\\
 2DSP\_10N2grp\_14 & \num{   0.0} & \num{   0.1} & \num{   0.0} & \num{         1} & \num{         1} & \num{         1} & \num{         1} & \num{         1} & \num{         1}\\
 2DSP\_10N2grp\_15 & \num{ 225.6} & \num{ 755.2} & \num{ 806.0} & \num{         1} & \num{         1} & \num{         1} & \num{        11} & \num{         7} & \num{         7}\\
 2DSP\_10N2grp\_16 & \num{  13.8} & \num{  21.9} & \num{  22.0} & \num{         1} & \num{         1} & \num{         1} & \num{        10} & \num{         7} & \num{         7}\\
 2DSP\_10N2grp\_17 & \num{   0.2} & \num{   0.2} & \num{   0.2} & \num{         1} & \num{         1} & \num{         1} & \num{         1} & \num{         1} & \num{         1}\\
 2DSP\_10N2grp\_19 & \num{   0.2} & \num{   0.2} & \num{   0.2} & \num{         1} & \num{         1} & \num{         1} & \num{         1} & \num{         1} & \num{         1}\\
  2DSP\_10N2grp\_2 & \num{  44.7} & \num{ 132.2} & \num{ 131.5} & \num{         1} & \num{         1} & \num{         1} & \num{         9} & \num{         6} & \num{         6}\\
 2DSP\_10N2grp\_21 & \num{   0.0} & \num{   0.0} & \num{   0.0} & \num{         1} & \num{         1} & \num{         1} & \num{         1} & \num{         1} & \num{         1}\\
 2DSP\_10N2grp\_22 & \num{   0.2} & \num{   0.2} & \num{   0.2} & \num{         1} & \num{         1} & \num{         1} & \num{         1} & \num{         1} & \num{         1}\\
 2DSP\_10N2grp\_23 & \num{  19.6} & \num{  46.7} & \num{  47.0} & \num{         4} & \num{         1} & \num{         1} & \num{        43} & \num{        15} & \num{        15}\\
 2DSP\_10N2grp\_24 & \num{   0.2} & \num{   0.2} & \num{   0.2} & \num{         1} & \num{         1} & \num{         1} & \num{         1} & \num{         1} & \num{         1}\\
 2DSP\_10N2grp\_25 & \num{   0.0} & \num{   0.0} & \num{   0.1} & \num{         1} & \num{         1} & \num{         1} & \num{         1} & \num{         1} & \num{         1}\\
 2DSP\_10N2grp\_26 & \num{   3.4} & \num{   9.3} & \num{   7.9} & \num{        23} & \num{        34} & \num{        19} & \num{        54} & \num{        32} & \num{        23}\\
 2DSP\_10N2grp\_27 & \num{   0.2} & \num{   0.2} & \num{   0.2} & \num{         1} & \num{         1} & \num{         1} & \num{         1} & \num{         1} & \num{         1}\\
 2DSP\_10N2grp\_28 & \num{ 440.9} & \num{2839.9} & \num{2851.9} & \num{         3} & \num{         3} & \num{         3} & \num{        28} & \num{        10} & \num{        10}\\
 2DSP\_10N2grp\_29 & \num{   7.5} & \num{ 195.1} & \num{ 194.0} & \num{         1} & \num{         1} & \num{         1} & \num{        10} & \num{        11} & \num{        11}\\
  2DSP\_10N2grp\_3 & \num{   0.9} & \num{   0.9} & \num{   1.0} & \num{         1} & \num{         1} & \num{         1} & \num{         1} & \num{         1} & \num{         1}\\
 2DSP\_10N2grp\_30 & \num{   0.8} & \num{   0.8} & \num{   0.8} & \num{         1} & \num{         1} & \num{         1} & \num{         1} & \num{         1} & \num{         1}\\
 2DSP\_10N2grp\_31 & \num{   0.1} & \num{   0.1} & \num{   0.1} & \num{         1} & \num{         1} & \num{         1} & \num{         1} & \num{         1} & \num{         1}\\
 2DSP\_10N2grp\_32 & \num{   3.6} & \num{ 147.0} & \num{ 147.0} & \num{         1} & \num{         1} & \num{         1} & \num{         8} & \num{         7} & \num{         7}\\
 2DSP\_10N2grp\_33 & \num{  28.3} & \num{ 147.1} & \num{ 147.0} & \num{         1} & \num{         2} & \num{         2} & \num{        10} & \num{        18} & \num{        18}\\
 2DSP\_10N2grp\_34 & \num{  70.1} & \num{  96.8} & \num{  96.7} & \num{         1} & \num{         1} & \num{         1} & \num{        17} & \num{         8} & \num{         8}\\
 2DSP\_10N2grp\_35 & \num{   1.1} & \num{   1.1} & \num{   1.1} & \num{         1} & \num{         1} & \num{         1} & \num{         1} & \num{         1} & \num{         1}\\
 2DSP\_10N2grp\_36 & \num{   2.2} & \num{  25.4} & \num{  25.5} & \num{         1} & \num{         1} & \num{         1} & \num{         8} & \num{        10} & \num{        10}\\
 2DSP\_10N2grp\_37 & \num{   4.7} & \num{  23.6} & \num{  22.3} & \num{        29} & \num{        32} & \num{        25} & \num{        93} & \num{        31} & \num{        29}\\
 2DSP\_10N2grp\_38 & \num{   0.2} & \num{   0.2} & \num{   0.2} & \num{         1} & \num{         1} & \num{         1} & \num{         1} & \num{         1} & \num{         1}\\
 2DSP\_10N2grp\_39 & \num{   8.8} & \num{  49.0} & \num{  34.8} & \num{        35} & \num{        52} & \num{        33} & \num{        96} & \num{        33} & \num{        27}\\
  2DSP\_10N2grp\_4 & \num{  21.7} & \num{  20.7} & \num{  20.7} & \num{         1} & \num{         2} & \num{         2} & \num{        18} & \num{        19} & \num{        19}\\
  2DSP\_10N2grp\_5 & \num{ 211.8} & \num{1683.5} & \num{1567.8} & \num{         4} & \num{        19} & \num{        18} & \num{        27} & \num{        73} & \num{        67}\\
  2DSP\_10N2grp\_6 & \num{   0.2} & \num{   0.2} & \num{   0.2} & \num{         1} & \num{         1} & \num{         1} & \num{         1} & \num{         1} & \num{         1}\\
  2DSP\_10N2grp\_8 & \num{   7.5} & \num{  21.0} & \num{  18.1} & \num{        27} & \num{        19} & \num{        17} & \num{        95} & \num{        14} & \num{        12}\\
  2DSP\_10N2grp\_9 & \num{ 205.9} & \num{2338.0} & \num{2553.4} & \num{         4} & \num{        13} & \num{        13} & \num{        37} & \num{        66} & \num{        66}\\
  2DSP\_15N3grp\_1 & \num{   1.9} & \num{   1.9} & \num{   1.9} & \num{         1} & \num{         1} & \num{         1} & \num{         1} & \num{         1} & \num{         1}\\
 2DSP\_15N3grp\_10 & \num{  78.6} & \num{ 482.2} & \num{ 357.6} & \num{        53} & \num{        85} & \num{        30} & \num{        99} & \num{        78} & \num{        59}\\
 2DSP\_15N3grp\_13 & \num{  85.8} & \num{1846.7} & \num{ 144.6} & \num{      1603} & \num{      1139} & \num{       229} & \num{      2395} & \num{       985} & \num{       104}\\
 2DSP\_15N3grp\_15 & \num{ 640.1} & \num{2375.4} & \num{ 851.4} & \num{       145} & \num{        40} & \num{        21} & \num{       339} & \num{        30} & \num{        16}\\
  2DSP\_15N3grp\_2 & \num{ 680.2} & \num{4333.1} & \num{1602.8} & \num{       199} & \num{       176} & \num{        60} & \num{       305} & \num{       122} & \num{        43}\\
 2DSP\_15N3grp\_21 & \num{   1.5} & \num{   1.5} & \num{   1.5} & \num{         1} & \num{         1} & \num{         1} & \num{         1} & \num{         1} & \num{         1}\\
 2DSP\_15N3grp\_22 & \num{   1.5} & \num{   1.5} & \num{   1.5} & \num{         1} & \num{         1} & \num{         1} & \num{         1} & \num{         1} & \num{         1}\\
 2DSP\_15N3grp\_26 & \num{1792.9} & \num{5489.1} & \num{5146.1} & \num{       107} & \num{       129} & \num{        66} & \num{       132} & \num{        72} & \num{        45}\\
 2DSP\_15N3grp\_29 & \num{ 290.4} & \num{ 562.8} & \num{ 404.4} & \num{        27} & \num{       137} & \num{        45} & \num{        63} & \num{       134} & \num{        41}\\
 2DSP\_15N3grp\_31 & \num{   1.4} & \num{   1.5} & \num{   1.5} & \num{         1} & \num{         1} & \num{         1} & \num{         1} & \num{         1} & \num{         1}\\
 2DSP\_15N3grp\_32 & \num{ 427.1} & \num{1541.7} & \num{ 541.7} & \num{      1270} & \num{      1041} & \num{       463} & \num{      1925} & \num{       752} & \num{       122}\\
 2DSP\_15N3grp\_33 & \num{ 120.0} & \num{ 604.6} & \num{ 213.4} & \num{       131} & \num{       137} & \num{        42} & \num{        76} & \num{       134} & \num{        50}\\
 2DSP\_15N3grp\_35 & \num{   1.4} & \num{   1.4} & \num{   1.4} & \num{         1} & \num{         1} & \num{         1} & \num{         1} & \num{         1} & \num{         1}\\
 2DSP\_15N3grp\_36 & \num{  21.4} & \num{ 126.6} & \num{  65.4} & \num{      1126} & \num{        52} & \num{        26} & \num{      1105} & \num{        87} & \num{        48}\\
 2DSP\_15N3grp\_37 & \num{  42.8} & \num{1215.0} & \num{ 835.2} & \num{       693} & \num{      1089} & \num{       404} & \num{      1410} & \num{       449} & \num{        76}\\
 2DSP\_15N3grp\_39 & \num{  42.8} & \num{1595.5} & \num{ 161.0} & \num{      2106} & \num{      5405} & \num{      2719} & \num{      4415} & \num{      1009} & \num{       110}\\
  2DSP\_15N3grp\_8 & \num{1267.2} & \num{4352.5} & \num{2097.0} & \num{       266} & \num{       176} & \num{       115} & \num{       324} & \num{        78} & \num{        40}\\
\midrule
mean (53 inst.) & \num{  12.6} & \num{  34.8} & \num{  27.7} & \num{      16.4} & \num{      16.5} & \num{      11.5} & \num{      15.9} & \num{      11.0} & \num{       8.1}\\
    \bottomrule
  \end{tabular*}
  \end{scriptsize}
\end{table}
\begin{table}[t]
  \caption{Detailed statistics for \scheduling instances that are solved by all Benders approaches.}
  \label{tab:fulloptsched}
  \centering
  \begin{scriptsize}
  \begin{tabular*}{\textwidth}{@{}l@{\;\;\extracolsep{\fill}}rrrrrrrrr@{}}
    \toprule
    & \multicolumn{3}{c}{running time} & \multicolumn{3}{c}{\#nodes} & \multicolumn{3}{c}{\#sepa}\\
    \cmidrule{2-4} \cmidrule{5-7} \cmidrule{8-10}
    instance & pool & EFrow & EFcons & pool & EFrow & EFcons & pool & EFrow & EFcons \\
    \midrule
 MSPS\_M3\_16J4\_0 & \num{ 208.1} & \num{  21.7} & \num{   5.6} & \num{    334963} & \num{      8098} & \num{      1537} & \num{    154053} & \num{      1279} & \num{        80}\\
 MSPS\_M3\_16J4\_1 & \num{ 745.8} & \num{ 469.6} & \num{  61.1} & \num{   1411173} & \num{    745902} & \num{     60954} & \num{    641905} & \num{     35281} & \num{       151}\\
MSPS\_M3\_16J4\_10 & \num{ 188.2} & \num{  20.9} & \num{   4.5} & \num{    251892} & \num{      4108} & \num{       791} & \num{    127872} & \num{       795} & \num{        49}\\
MSPS\_M3\_16J4\_11 & \num{ 744.0} & \num{ 153.8} & \num{  16.2} & \num{   1498201} & \num{    176372} & \num{      6418} & \num{    658155} & \num{     23631} & \num{        70}\\
MSPS\_M3\_16J4\_12 & \num{1430.1} & \num{  60.3} & \num{   9.7} & \num{   2589118} & \num{     55353} & \num{      1912} & \num{   1305023} & \num{     21638} & \num{       183}\\
MSPS\_M3\_16J4\_13 & \num{ 496.1} & \num{  48.4} & \num{  11.3} & \num{    888909} & \num{     42583} & \num{      2131} & \num{    380819} & \num{     16082} & \num{        51}\\
MSPS\_M3\_16J4\_14 & \num{2366.7} & \num{ 134.8} & \num{  26.1} & \num{   4102236} & \num{    168881} & \num{      4468} & \num{   1674821} & \num{     56647} & \num{       143}\\
MSPS\_M3\_16J4\_15 & \num{ 700.9} & \num{  64.4} & \num{  18.5} & \num{   1399138} & \num{     80900} & \num{      4990} & \num{    727886} & \num{     29817} & \num{        75}\\
MSPS\_M3\_16J4\_16 & \num{ 152.8} & \num{  15.9} & \num{   2.1} & \num{    225335} & \num{      9492} & \num{       444} & \num{    100422} & \num{      3857} & \num{        89}\\
MSPS\_M3\_16J4\_17 & \num{2961.8} & \num{ 127.8} & \num{  19.5} & \num{   5085046} & \num{    140884} & \num{      3072} & \num{   2564053} & \num{     63479} & \num{        80}\\
MSPS\_M3\_16J4\_18 & \num{1772.2} & \num{ 312.1} & \num{  11.6} & \num{   2972314} & \num{    493631} & \num{      2410} & \num{   1478175} & \num{    242596} & \num{       294}\\
MSPS\_M3\_16J4\_19 & \num{2382.8} & \num{ 241.4} & \num{  24.0} & \num{   4312101} & \num{    326184} & \num{      5301} & \num{   2146882} & \num{    138837} & \num{       293}\\
 MSPS\_M3\_16J4\_2 & \num{2367.6} & \num{ 495.0} & \num{  12.7} & \num{   3962974} & \num{    688217} & \num{      2092} & \num{   1934204} & \num{    312764} & \num{       222}\\
MSPS\_M3\_16J4\_20 & \num{ 715.1} & \num{  20.4} & \num{   4.5} & \num{   1302010} & \num{      8942} & \num{      1231} & \num{    618990} & \num{      2650} & \num{        58}\\
MSPS\_M3\_16J4\_21 & \num{ 635.2} & \num{  36.9} & \num{  11.0} & \num{   1066923} & \num{     30551} & \num{      1595} & \num{    500595} & \num{      6964} & \num{       107}\\
MSPS\_M3\_16J4\_22 & \num{ 960.4} & \num{ 100.4} & \num{  20.5} & \num{   1883006} & \num{    115854} & \num{      4575} & \num{    942525} & \num{     41231} & \num{       127}\\
MSPS\_M3\_16J4\_23 & \num{ 503.5} & \num{  36.5} & \num{   8.9} & \num{    991814} & \num{     21650} & \num{      2049} & \num{    492281} & \num{      6835} & \num{        56}\\
MSPS\_M3\_16J4\_24 & \num{1031.4} & \num{  38.0} & \num{  14.3} & \num{   1833236} & \num{     24079} & \num{      2076} & \num{    906884} & \num{      8125} & \num{       100}\\
MSPS\_M3\_16J4\_25 & \num{ 620.2} & \num{  28.9} & \num{   2.7} & \num{   1053714} & \num{     11180} & \num{       724} & \num{    492034} & \num{      3310} & \num{        50}\\
MSPS\_M3\_16J4\_26 & \num{ 761.2} & \num{  42.1} & \num{  15.9} & \num{   1319166} & \num{     15211} & \num{      3147} & \num{    688770} & \num{      2945} & \num{        76}\\
MSPS\_M3\_16J4\_27 & \num{ 561.8} & \num{  45.2} & \num{  12.2} & \num{   1153721} & \num{     39555} & \num{      2459} & \num{    590868} & \num{      1918} & \num{        53}\\
MSPS\_M3\_16J4\_28 & \num{ 872.2} & \num{  72.7} & \num{  25.0} & \num{   1813662} & \num{     62549} & \num{     11210} & \num{    743536} & \num{     16238} & \num{        46}\\
MSPS\_M3\_16J4\_29 & \num{ 232.7} & \num{  16.9} & \num{   2.1} & \num{    382419} & \num{      4349} & \num{       303} & \num{    183254} & \num{       907} & \num{        56}\\
 MSPS\_M3\_16J4\_3 & \num{1513.7} & \num{  84.0} & \num{  13.4} & \num{   2736777} & \num{     57898} & \num{      2261} & \num{   1235881} & \num{      6475} & \num{        80}\\
MSPS\_M3\_16J4\_30 & \num{1111.2} & \num{ 160.7} & \num{  15.5} & \num{   2235635} & \num{    232802} & \num{      4278} & \num{   1104867} & \num{     97622} & \num{       294}\\
MSPS\_M3\_16J4\_31 & \num{ 131.4} & \num{  24.2} & \num{   2.9} & \num{    198873} & \num{     15556} & \num{       880} & \num{     82820} & \num{      1593} & \num{        42}\\
MSPS\_M3\_16J4\_32 & \num{ 230.8} & \num{  19.5} & \num{   4.6} & \num{    341899} & \num{      7219} & \num{       956} & \num{    177159} & \num{      2148} & \num{        58}\\
MSPS\_M3\_16J4\_33 & \num{ 614.9} & \num{  49.4} & \num{   9.2} & \num{   1179530} & \num{     34978} & \num{      1818} & \num{    579800} & \num{      4883} & \num{        91}\\
MSPS\_M3\_16J4\_34 & \num{ 629.9} & \num{  22.4} & \num{   7.3} & \num{   1030701} & \num{      7275} & \num{       902} & \num{    488129} & \num{      2561} & \num{       219}\\
MSPS\_M3\_16J4\_35 & \num{ 951.8} & \num{ 341.1} & \num{  23.9} & \num{   1819960} & \num{    520265} & \num{      6926} & \num{    841374} & \num{    117639} & \num{        72}\\
MSPS\_M3\_16J4\_36 & \num{ 905.1} & \num{  41.8} & \num{   8.0} & \num{   1571586} & \num{     31310} & \num{      1702} & \num{    702296} & \num{      6251} & \num{       168}\\
MSPS\_M3\_16J4\_37 & \num{ 349.5} & \num{  57.2} & \num{   9.5} & \num{    622188} & \num{     63987} & \num{      1958} & \num{    318018} & \num{     22120} & \num{        54}\\
MSPS\_M3\_16J4\_38 & \num{2069.0} & \num{ 152.6} & \num{  15.0} & \num{   3325297} & \num{    179661} & \num{      2816} & \num{   1631977} & \num{     78653} & \num{       161}\\
MSPS\_M3\_16J4\_39 & \num{  93.6} & \num{  17.0} & \num{   2.7} & \num{    136274} & \num{      8520} & \num{       385} & \num{     58457} & \num{      3638} & \num{        67}\\
 MSPS\_M3\_16J4\_4 & \num{1348.2} & \num{  37.1} & \num{  14.5} & \num{   2553764} & \num{     21753} & \num{      3077} & \num{   1250382} & \num{      6591} & \num{       119}\\
 MSPS\_M3\_16J4\_5 & \num{1436.2} & \num{ 488.1} & \num{  30.2} & \num{   2852482} & \num{    662332} & \num{      9764} & \num{   1361762} & \num{    287749} & \num{       163}\\
 MSPS\_M3\_16J4\_6 & \num{1410.5} & \num{  32.1} & \num{   9.5} & \num{   2212063} & \num{     16124} & \num{      1488} & \num{   1101449} & \num{      6457} & \num{        90}\\
 MSPS\_M3\_16J4\_7 & \num{ 776.2} & \num{ 283.2} & \num{  31.6} & \num{   1477659} & \num{    410993} & \num{     14588} & \num{    700793} & \num{     54494} & \num{        98}\\
 MSPS\_M3\_16J4\_8 & \num{1241.8} & \num{ 178.3} & \num{  20.9} & \num{   2525524} & \num{    254203} & \num{      7324} & \num{   1244443} & \num{     93452} & \num{       135}\\
 MSPS\_M3\_16J4\_9 & \num{1230.6} & \num{ 337.6} & \num{  20.9} & \num{   2462718} & \num{    514077} & \num{      4626} & \num{   1203918} & \num{    203483} & \num{        75}\\
MSPS\_M3\_19J4\_10 & \num{4482.8} & \num{ 283.2} & \num{  14.2} & \num{   5458847} & \num{    284089} & \num{      5499} & \num{   2021252} & \num{     32560} & \num{       109}\\
MSPS\_M3\_19J4\_16 & \num{6113.0} & \num{  24.9} & \num{   8.8} & \num{   6813189} & \num{      8062} & \num{      1084} & \num{   2820143} & \num{      1536} & \num{        59}\\
MSPS\_M3\_19J4\_29 & \num{5236.9} & \num{  45.2} & \num{  14.5} & \num{   6664102} & \num{     20042} & \num{      2634} & \num{   2986595} & \num{      4093} & \num{        87}\\
MSPS\_M3\_19J4\_31 & \num{3745.1} & \num{ 150.9} & \num{  74.4} & \num{   5047269} & \num{    120995} & \num{     34202} & \num{   1805159} & \num{      5402} & \num{       140}\\
MSPS\_M3\_19J4\_32 & \num{5506.0} & \num{ 135.0} & \num{  15.0} & \num{   6776613} & \num{    101483} & \num{      2179} & \num{   3215502} & \num{     36872} & \num{        90}\\
\midrule
mean (45 inst.) & \num{ 912.3} & \num{  72.7} & \num{  12.0} & \num{ 1550407.5} & \num{   55351.8} & \num{    2649.2} & \num{  726246.4} & \num{   13587.5} & \num{      95.6}\\
    \bottomrule
  \end{tabular*}
  \end{scriptsize}
\end{table}

\end{document}